

\input amstex
\documentstyle{amsppt}

\NoRunningHeads
\magnification=1100
\baselineskip=29pt
\parskip 9pt
\pagewidth{5.2in}
\pageheight{7.2in}
\TagsOnRight
\CenteredTagsOnSplits

\def\aut{\#_{\text{aut}}}
\def\Aut{\operatorname{Aut}}
\def\alp{\alpha}

\def\Bl{\bigl(}
\def\Br{\bigr)}
\def\BBl{\Bigl(}
\def\BBr{\Bigr)}

\def\codim{\operatorname{codim}}
\def\cdott{\!\cdot\!}
\def\ctil{\tilde{C}}
\def\CC{{\Bbb C}}

\def\Del{\Delta}
\def\dhct{D^h_{\ctil}}
\def\dh{D^h}
\def\dual{^{\vee}}

\def\endo{\Cal{E}nd^0}

\let\eps=\varepsilon
\def\E{\Cal{E}}

\def\Ext{\operatorname{Ext}}

\def\eo{E^{(1)}}
\def\et{E^{(2)}}
\def\ez{E^{(0)}}
\def\erdi{\eta(r,d,I)}

\def\F{\Cal{F}}

\def\Gam{{\Upsilon}}
\def\Gamlam{\Upsilon_{\Lambda}}

\def\Hom{\text{Hom}}
\def\half{{1\over 2}}
\def\hdc{|hH_C(-h\Del)|}
\def\he{H(\eps)}
\def\hed{H(\eps)_{|\Del}}

\let\lra=\longrightarrow
\let\len=\ell
\def\LL{\Cal{L}}
\def\llam{_{\Lambda}}

\def\mh{\!:\!}

\def\mde{\MM^{d,\eps}}
\def\mod{\#_{\text{mod}}}
\def\modl{\#_{\text{mod}}^{\text{loc}}}
\def\mrdx{\MM^{r,d}_X}
\let\mxrd=\mrdx
\def\Mx{\MM^{r,d}_X(I,H)}
\def\mtdx{\MM^{2,d}_X}

\def\MM{{\frak M}}
\def\Md{\frak M^d}
\let\md=\Md
\def\mapright#1{\,\smash{\mathop{\lra}\limits^{#1}}\,}

\def\nue{\sim}
\def\nde{{\frak N}^{d,\eps}}

\def\OO{\Cal{O}}

\def\otimess{\!\otimes\!}

\def\oplu#1#2{\smash{\mathop{\oplus}\limits^{#1}_{#2}}}
\def\ome{\omega}

\def\pl{\!+\!}
\def\plu{\pl}
\def\PP{\bold{P}}
\def\Pic{\operatorname{Pic}}
\def\pri{^{\prime}}
\def\Px{P_{\xi}}
\def\ps{p\sta}
\def\pxs{p_X\sta}

\def\pr{\text{pr}}
\def\proof{{\it Proof.\ }}
\let\Proof=\proof

\let\qqed=\qedd
\def\endpf{\hfill\qed\vskip8pt}

\def\Quv{\bold{Quot}^c_V}
\def\QQ{{\Bbb Q}}

\def\fretl{\Upsilon_{\Lambda}}

\def\rank{\operatorname{rk}}
\def\rk{\operatorname{rk}}
\def\R{\Cal{R}}
\def\Rr{\frak R^{r,d}_{e,I}}

\def\rd{^{r,d}_{e,I}}
\def\rdi{^{r_i,d_i}_{e,I_i}}
\def\rdz{^{r,d}_{\mu,I}}
\def\rdp{^{r,d\pri}_{e,I}}
\def\rdo{^{r,d}_{e_1,I}}
\def\rdt{^{r,d}_{e_2,I}}
\def\rdm{^{r,d}_{\mu,I}}
\def\RR{{\Bbb R}}

\def\sslo{\Cal{S}_{\Lambda,\ome}}
\def\sing{\operatorname{Sing}}
\def\spec{\operatorname{Spec\,}}
\let\sub=\subseteq
\def\Sig{\Sigma}
\def\sigm{\Sigma^{-}}
\def\Sigp{\Sigma^+}
\def\sigp{\Sigma^{+}}
\def\sta{^{\ast}}
\def\SL{\Cal{SL}}
\def\she{\Cal{S}}
\def\supp{\operatorname{supp}}
\def\SS{\Cal{S}}

\def\te{\tilde E}
\def\timess{\!\times\!}
\def\timec{\!\times_C\!}
\def\thh{\tilde{\theta}^h}
\def\tilL{\tilde L}

\def\U{\Cal{U}}
\def\uplu#1{^{\oplus #1}}

\def\upr{^{\oplus r}}
\def\upmo{^{-1}}

\def\vec{\frak A}
\def\vecx{\Cal{V}}

\def\Wr{\frak A^{r,d}_{e,I}}
\def\Wrp{\frak A^{r,d\pri}_{e\pri,I\pri}}
\def\Wrt{\frak A^{r,d_2}_{e_2,I_2}}

\def\xx{\chi}
\def\Xix{\Xi_1({\frak x})}

\def\ZZ{{\Bbb Z}}

\let\pro=\proclaim
\let\endpro=\endproclaim

\topmatter
\title
Moduli of high rank Vector bundles over Surfaces
\endtitle

\author
David Gieseker \\
Jun Li
\endauthor
\thanks
\endthanks

\thanks The first author was partially supported
by NSF grant DMS-9305657 and the second author was
partially supported by NSF grant DMS-9307892
\endthanks

\affil
Mathematics Department\\
University of California, Los Angeles
\endaffil

\address
Mathematics Department\\
University of California, Los Angeles\\
Los Angeles, CA 90024\\
\endaddress

\email
dag@math.ucla.edu//jli@math.ucla.edu
\endemail

\endtopmatter
\document

\head 0. Introduction
\endhead

The purpose of this work
is to apply the degeneration theory developed in [GL] to study the
moduli space of stable vector bundles of arbitrary rank on any smooth
algebraic surface (over $\CC$). We will show that most of the recent
progress in understanding
moduli of rank two vector bundles can be carried over
to high rank cases.

After introducing the notion of
stable vector bundles, the first author constructed the moduli
schemes of vector bundles on surfaces. He showed
that for any smooth algebraic surface $X$\ with ample divisor
$H$\ and line bundle $I$\ on $X$,
there is a coarse moduli scheme $\mrdx(I,H)$\ parameterizing (modulo
equivalence relation) the set of all $H$-semistable rank $r$\ torsion free
sheaves $E$\ on $X$\ with $\det E=I$\ and $c_2(E)=d$.
Since then, many mathematicians have studied the geometry of this moduli
space, especially for rank two case. To cite a few,
Maruyama, Taubes and
the first author showed that the moduli space $\mtdx$\ ($=\mrdx(I,H)$) is
non-empty when $d$\ is large.
Moduli spaces of vector bundles of some special surfaces have been
studied also.

The deep understanding of $\mrdx$\ for arbitrary $X$\ and $r=2$\
begins with Donaldson's generic smoothness result. Roughly speaking,
Donaldson [Do], (later generalized by Friedman [Fr] and K. Zhu [Zh])
showed that when $d$\ is large enough, then the singular locus $\sing\Bl
\mtdx\Br$\ of $\mtdx$\ is a proper subset of $\mtdx$\
and its codimension in $\mtdx$\ increases linearly
in $d$. This theorem indicates that the moduli
$\mtdx$\ behaves as expected when the second Chern class $d$\ is
large. Later, using general deformation theory, the second author proved
that $\mtdx$\ is normal, and has local complete intersection (l.c.i.)
singularities at stable
sheaves provided $d$\ is large [L2]. He also showed that when $X$\ is a
surface of general type satisfying some mild technical conditions, then
$\mtdx$\ is of general type for $d\gg 0$\ [L2]. In our paper [GL],
we also proved that $\mtdx$\ is irreducible if $d$\ is large.

In this and subsequent papers, we shall show that
the geometry
of $\mtdx$\ and the geometry of $\mrdx$, $r\geq 3$, is rather
similar. The main obstacle in doing so is the lack of
an analogy of the generic
smoothness result in high rank case.
In this paper,
we will use the degeneration of moduli developed in [GL] to establish the
following main technical theorem:

\pro{Theorem 0.1}{
Let $X$\ be a smooth algebraic surface, $H$\ an ample line bundle and
$I$\ a line bundle on $X$. Let $r\geq 2$\ be any integer. Then for any
constant $C_1$\ and any divisor $D\sub X$, there is an $N$\ such that
whenever $d\geq N$, then we have
$$\dim\bigl\{E\in\mrdx\mid \Ext^0(E,E(D))^0\ne\{0\}\bigr\}
\leq\eta(r,d,I)-C_1,
$$
where $\eta(r,d,I)=2rd-(r-1)I^2-(r^2-1)\xx(\OO_X)$\ is the expected
dimension of $\mrdx$\ $(=\mrdx(I,H))$\ and the superscript $0$\
stands for the traceless part of $\Ext^i(\cdot,\cdot)$.
\endpro

According to [At][Mu], $\mrdx$\ is regular at $E$\ if $E$\ is stable
and $\Ext^2(E,E)^0=\{0\}$. As to the subset
of strictly semistable sheaves in $\mrdx$, it is easy to show that
its dimension is much less
than $\eta(r,d,I)-C_1$\ when $d$\ is large. After applying theorem 0.1
to the divisor $D=K_X$\ and using the Serre duality, we
conclude that for $d$\ sufficiently large,
$$\dim\sing\Bl\mrdx\Br\leq\eta(r,d,I)-C_1.
$$
On the other hand, based on deformation theory, each component of
$\mrdx$\ has dimension at least $\eta(r,d,I)$. Thus, we have
proved the following theorem.

\pro{Theorem 0.2}
Let $X$\ be a smooth algebraic surface, $H$\ an ample line bundle and
$I$\ a line bundle on $X$. Let $r\geq 2$\ be any integer. Then for any
constant $C_1$, there is an $N$\ such that whenever $d\geq N$, then
$\mrdx$\ has pure dimension $\erdi$\ and further,
$$\codim\Bl\sing\Bl\mrdx\Br,\mrdx\Br\geq C_1.
$$
\endpro

Once we have settled the generic smoothness result, we can generalize some
other properties of $\mtdx$\ to high rank case. In this paper, we will
prove

\pro{Theorem 0.3}
With the notation as in theorem 0.2, then there is an $N$\ such that
whenever $d\geq N$, then
\roster
\item $\mrdx$\ is normal. Further, if $s\in\mrdx$\
is a closed point
corresponding to a stable sheaf, then $\mrdx$\ is a local
complete intersection at $s$;
\item The set of locally free $\mu$-stable sheaves $\Bl\mrdx\Br^{\mu}\sub
\mrdx$\ is dense in $\mrdx$\ and
\item For any polarizations $H_1$\ and $H_2$\ of $X$, the moduli
$\mrdx(I,H_1)$\ is birational to $\mrdx(I,H_2)$. \text (In this case, $N$\
depends on both $H_1$\ and $H_2$.\text )
\endroster
\endpro

To illustrate the idea of the proof of our main theorem (theorem 0.1), let us
first recall the degeneration of moduli $\mrdx$\ constructed in [GL].
Let $0\in C\sub\spec\CC[t]$\ be a smooth curve that functions as a parameter
space and let $Z\to C$\ be a family of surfaces that is the result of
blowing-up $X\times C$\ along $\Sigma\times\{0\}$, where $\Sigma\in|H|$\
is a smooth very ample divisor. Clearly, $Z_t=\pi^{-1}(t)$\ is $X$\
and $Z_0=X\cup \Del$, where $\Del$\ is a ruled surface over $\Sigma$. Over
$C\sta=C\setminus\{0\}$, we have a constant family $\mrdx\times C\sta$.
In [GL], we have constructed completions of $\mrdx\times C\sta$\
over $C$. These completions depend on the choice of ample
divisors on $Z$. The ample divisor which we will use is a multiple of
the $\QQ$-divisor $p_X\sta H(-(1-\eps)\Del)$\ that depends on the rational
$\eps\in (0,{1\over 2})$. We denote this completion by
$\mde$. There is a nice description of closed points of the special
fiber $\mde_0$: Any point of $\mde_0$\ corresponds uniquely to an
equivalence class of semistable sheaves on $Z_0$.

Now let $D\sub X$\ be any divisor and ${\Cal N}\sub\mrdx$\
be the set of sheaves $E$\ such that
$$\Hom(E,E(D))^0\ne\{0\}.\tag 0.1
$$
Put $\nde\sub\mde$\ be the closure of ${\Cal N}
\times C\sta$\ in $\mde$. To show that for any constant $C_1$\
and large $d$\ we have
$$\dim{\Cal N}\leq\eta(r,d,I)-C_1,
$$
it suffices to show
$$\dim\nde_0\leq\eta(r,d,I)-C_1.\tag 0.2
$$
Now let $E\in\nde_0$\ be any sheaf. Note that $E$\ is a limit of sheaves
in $\Cal N$\ and that sheaves in $\Cal N$\ satisfy (0.1). So
by semicontinuity theorem, for any invertible sheaf ${\Cal L}$\
on $Z$\ such that ${\Cal L}_{|Z_t}\cong \OO_X(D)$, we have
$$\Hom_{Z_0}(E,E\otimes{\Cal L}_{Z_0})^0\ne\{0\}.
$$
In particular,
if we choose ${\Cal L}$\ to be $\pxs\OO_X(D)(-k\Del)$,
where $p_X\mh Z\to X$\ is the projection, we get
$$\Ext^0_{Z_0}(E,E\otimes\pxs\OO_X(D)(-k\Del))^0\ne\{0\},\quad
\forall k\in\ZZ,\ E\in\nde_0. \tag 0.3
$$
Since $E$\ is semistable, $E_{|X}$\ and $E_{|\Del}$\ as sheaves on $X$\
and $\Del$\ respectively will satisfy some weak stability conditions.
(For simplicity, here we assume $E$\ is locally free.) On the other hand,
for large $k$, the non-vanishing of
$$\Ext^0_X(E_{|X},E_{|X}(D-k\Sigma))^0\tag 0.4
$$
will force $E_{|X}$\ to be very unstable. Therefore, we
can choose a $k>0$\ (independent of $d$, $\eps$\
and $\Cal N$) such that (0.4) is always trivial.
Thus (0.3) will force
$$\Ext^0_{\Del}(E_{|\Del},E_{|\Del}\otimes\pxs\OO_X(D)(k\Sigma))^0\ne\{0\}.
\tag 0.5
$$
(0.5) certainly is possible for sheaves over $\Del$. However, if
we can show that the number of moduli of the set of sheaves
$F$\ (over $\Del$) satisfying (0.5) is strictly less than
$$\hbox{the number of moduli of}\ \{E_{\Del}\mid E_\in \mde_0\}-C_1,
$$
then $\codim(\nde_0,\mde_0)\geq C_1$, which is exactly what we need.
Therefore, the proof of theorem 0.1 is reduced to the proof of the
following theorem.

\pro{Theorem 0.4}
Let $X$\ be any ruled surface and let $H$\ and $I$\ be as in theorem 0.1. Then
for any integer $r$, any divisor $D\sub X$\ and any constant $C_1$, there is a
constant $N$\ such that for $d\geq N$,
$$\dim\bigl\{ E\in\mxrd\mid \Ext^0(E,E(D))^0\ne\{0\}\bigr\}
\leq\eta(r,d,I)-C_1.
$$
\endpro

The advantage of working with a ruled surface lies in the fact that
every vector bundle on ruled surface can be constructed
explicitly as follows: Let $X=
\Del$\ and let $E$\ be a vector bundle on $\Del$. For simplicity, we assume
for general fiber $P_{\xi}$\ of
$\pi\mh\Del\to\Sigma$, the restriction sheaf
$E_{|P_{\xi}}\cong\OO_{P_{\xi}}^{\oplus r}$.
 Then there is a unique rank $r$\ vector bundle $V$\
on $\Sigma$\ and a sheaf $F$\ supported on a finite number of
fibers of $\pi$\ such that
$$0\lra E\lra \pi\sta V\mapright{\varphi} F\lra 0
$$
is exact. When $E$\ is general, $F$\ is of the form $\oplus \OO_{P_i}
(1)$, where $P_i$\ are fibers of $\pi$. Thus
the condition under which $E$\ admits traceless homomorphism
$E\to E(D)$\ can be interpreted in terms of the location of $P_i$'s
and the choice of homomorphism $\varphi$. The argument to carry
out this approach is rather straightforward though quite technical
and will occupy the first
section of this paper. In \S 2, we will review the degeneration construction
and use it to prove theorem 0.1. The theorems 0.2-0.4 will be proved
in \S 3. We remark that after the completion of the initial version of this
work, O'Grady has improved our results in his paper [OG].

\subhead  Conventions and preliminaries
\endsubhead

All schemes are defined over the field of complex number $\CC$\
and are of finite type. All points are closed points unless otherwise
is mentioned. We shall always identify a vector bundle with its
sheaf of sections. If $I$\ and $J$\ are two line bundles on surface, then
we denote by $I\cdot J$\ the intersection $c_1(I)\cdot c_1(J)$\ and
$I^2$\ the self-intersection $c_1(I)\cdot c_1(I)$.
We will use $\nue$\ to denote the numerical equivalence of divisors
(line bundles). For coherent sheaf $F$, we denote by $\rk(F)$\ the
rank of $F$. In case $F$\ is supported on finite number of points on $X$,
we denote by $\len(F)$\ the length of $F$. If $p$\ and $q$\
are two polynomials with real coefficients, we say $p\succ q$\
(resp. $p\succeq q$) if $p(n)>q(n)$\ (resp. $p(n)\geq q(n)$)
for all $n\gg0$.

In the following, $X$\ will always denote a smooth projective surface. Let
$H$\ be a very ample line bundle on $X$. For any sheaf $E$\ on $X$, we denote
by $\xx_E$\ the Poincare polynomial of $E$, namely, $\xx_E(n)=\xx(E(n))$,
$E(n)=E\otimess H^{\otimes n}$, and denote by $p_E$\ the polynomial
${1\over \rk(E)}\xx_E$\ when $\rk(E)\ne0$. Unless the contrary is mentioned,
the degree of a sheaf $E$\ is $c_1(E)\cdot H$. We recall the notion of
stability:

\pro{Definition 0.4}
A sheaf $E$\ on $X$\ is said to be stable (resp. semistable) with respect to
$H$\ if $E$\ is coherent, torsion free and if one of the following two
equivalent conditions hold:\hfil\break
1. Whenever $F\subset E$\ is a proper subsheaf, then
$p_F\prec p_E\quad (\text{resp.}\ p_F\preceq p_E)$;\hfil\break
2. Whenever $E\to Q$\ is a quotient sheaf, $\rk(Q)>0$, then
$p_E\prec p_Q\quad (\text{resp.}\ p_E\preceq p_Q)$.
}
When $E$\ is a torsion free coherent sheaf on $X$, we define the slope
$\mu(E)={1\over \rk(E)}\deg E$.
\endpro

\pro{Definition 0.5}
Let $e$\ be a constant. The sheaf $E$\ is said to be $e$-stable if one
of the following two equivalent conditions hold:
\hfil\break
1. Whenever $F\subset E$\ is a subsheaf with $0<\rk(F)<\rk(E)$, then
$\mu(F)<\mu(E)+{1\over \rk (F)}\sqrt{H^2}\cdot e$;\hfil\break
2. Whenever $E\to Q$\ is a quotient sheaf with $0<\rk(Q)<\rk(E)$, then
$\mu(E)<\mu(Q)+{1\over \rk(Q)}\sqrt{H^2}\cdot e$.\hfil\break
We call $E$\ $\mu$-stable if $E$\ is $e$-stable with $e=0$. When the
strict inequality are replaced by $\leq$, then we call $E$\ $e$-semistable.
\endpro

Let $W\to S$\ be a flat morphism and let $E\to W$\ be any sheaf on $W$. For
any closed $s\in S$, we will use $W_s$\ to denote the fiber of $W$\ over
$s$\ and use $E_s$\ to denote the restriction of $E$\ to $W_s$. For any
subscheme $T\sub W$, we denote by $E_{|T}$\ the restriction of $E$\ to $T$.
We shall adopt the following convention:
If $\R$\ is a set of sheaves on $X$, then the number of moduli of
$\R$\ is the smallest integer $m$\ so that there are countably many
schemes (of finite types) of dimension at most $m$, say
$S_1, S_2, \cdots$, and flat family of sheaves $E_{S_1}, E_{S_2}, \cdots$\
on $X\timess S_1, X\timess S_2,\cdots$\ respectively of which the following
holds: For any $F\in \R$, there is a closed $s\in S_k$\ for some $k$\
such that $ F\cong E_{S_k, s}$.
We will denote by $\mod(\R)$\ the number of moduli of $\R$. In case $R$\ is
a scheme parameterizing a family of sheaves and $t\in R$, then we denote
by $\modl(R, [t])$\ the number of moduli of sheaves parameterized by
the germ of $R$\ at $r$. In particular, we write $\modl(E)$, where $E$\
is any sheaf, for $\modl({\Cal Q}, [E])$\ ($=$\ the number of moduli of
the set of all ``small'' deformations of $E$), where ${\Cal Q}$\ is a
Grothendieck's Quot-scheme [Gr] that contains all deformations
of $E$\ as quotient sheaves of some appropriate locally free sheaf. Another
notion we use frequently is $\aut(E)=\dim \Aut(E)$, where $\Aut(E)$\ is the
group of automorphism of $E$. Note $\aut(E)=h^0({\Cal E}nd(E))$. When $\R$\ is
a
set of sheaves, then $\aut(\R)=\max\{\aut(E)\mid E\in \R\}$.

\head 1. Vector bundles on ruled surface
\endhead

\def\vec{\frak A}
\def\Rr{\vec\rd}
\def\Sigp{\Sigma^+}

The purpose of this section is to prove an analogy of theorem
0.1 for ruled surface $\Del$. Before giving the precise statement
of the theorem, we first introduce some notation. Let $\Sig$\ be
a smooth curve and let $\pi\mh \Del\to \Sig$\ be a ruled surface.
For simplicity, we assume $\Del$\ is the projective bundle of a
direct sum of a trivial line bundle with a very ample line bundle
(over $\Sig$). Hence $\pi\mh \Del\to\Sig$\
has a unique section $\sigm$\ with $\sigm\cdott\sigm<0$\ and has many sections
with positive self-intersection. We choose one such section
and denote it by $\Sigp$. By assumption,
$|\sigp|$\ is base point free. Let $H$\ be an ample line bundle on $\Del$\
that is numerically equivalent to (denoted by $\sim$) $a\sigp+bP_{\xi}$,
where $P_{\xi}$\ is a general fiber of $\pi$.
Let $e$\ be a constant, let $I$\ be a line bundle on $\Del$\
and let $D$\ be any divisor on $\Del$. In this section, we will study the
set $\vec^{r,d}_{e,I,H}$\ of all $e$-semistable (with respect
to $H$) rank $r$\ locally free sheaves $E$\ with $\det E=I$\ and $c_2(E)=d$\
and the set
$$\vec^{r,d}_{e,I,H}(D)=\{E\in \vec^{r,d}_{e,I,H}\mid
\Hom(E,E(D))^0\ne\{0\}\}.
$$
Here and in the following, the superscript 0 always stands for the
traceless part of the group or sheaf.
For technical reasons, we will choose $H$\ to be very close to $\Sigp$\
in the sense that $b/a$\ is very small. With the
choice of $H$\ understood, we will not build $H$\ into the notation and will
write $\vec\rd$\ (resp. $\vec\rd(D)$)
for $\vec^{r,d}_{e,I,H}$\ (resp. $\vec^{r,d}_{e,I,H}(D)$).
We will also use
$\eta_{\Del}(E)=\eta_{\Del}(\rank(E),c_2(E),c_1(E))$\ to denote
the number
$$\eta_{\Del}(r,d,I)=2rd-(r-1) I^2-(r^2-1)\xx(\OO_{\Del}).\tag 1.1
$$
$\eta_{\Del}(r,d,I)$\ is the expected dimension of $\vec\rd$.
Because in this section we work solely with the surface $\Del$, we will
simply write $\eta$\ for $\eta_{\Del}$.
The theorem we will prove in this section is the following.

\pro{Theorem 1.1}
Given $r$\ and $\Del$, there is an $\eps_0>0$\ depending on $r$\
and $\Del$\ of which the following holds: For any ample
divisor $H\nue a\sigp+b\Px$\ with $b/a<\eps_0$\ and for any choice of
constants $e$, $C$\ and divisor $D\sub \Del$, $I\in \Pic(\Del)$,
there is an integer $N$\ such that whenever $d\geq N$, then we have
$$\mod\Bl\vec\rd(D)\Br\leq \eta(r,d,I)-C.\tag 1.2
$$
\endpro

The advantage of working with ruled surface lies on having a powerful
structure theorem of torsion free sheaves on $\Del$.
Let $E$\ be any torsion free sheaf of rank $r$.
By Grothendieck's splitting theorem, its restriction to a generic
fiber $P_{\xi}$\ has the form
$$E_{|P_{\xi}}\cong\oplu{n}{i=1}\OO_{P_{\xi}}(\alpha_i)^{\oplus r_i},
\quad \alpha_1>\cdots>\alpha_n.\tag 1.3
$$
In the following, we
call $\alp=(\alpha_1^{\oplus r_1},\cdots,\alpha_n^{\oplus r_n})$\
the generic fiber type of $E$.
(The integer sequence $\{\alpha_i\}$\ is always assumed to be strictly
decreasing.) We let $\len(\alpha)=\sum^n_{i=1} r_i\alpha_i$.
Clearly, $r=\sum^n_{i=1}r_i$\ and
further, when $\det E=I$\ and $\deg I_{|P_{\xi}}=m$, then
$m=\len(\alpha)$.
$\Rr$\ can be divided into strata according to the generic
fiber types of individual vector bundles.
Let $r\in {\Bbb N}$\ and $I\in \Pic(\Del)$\ be fixed.
Without loss of generality, we can assume $0\leq \deg I_{|P_{\xi}}\leq
r-1$. Let $m=\deg I_{|P_{\xi}}$\ and let $1_m$\ be the fiber type
$(1^{\oplus m},0^{\oplus(r-m)})$. For any fiber type $\alpha$\
with $\len(\alpha)=m$, we let
$$\Rr(\alp)=\{E\in\Rr\mid E\ \hbox{has generic fiber type}\ \alp\}.
$$
The first observation we have is that except for $\alpha=1_m$, none
of $\mod\bigl(\Rr(\alpha)\bigr)$\ are close to
$\eta(r,d,I)$. More precisely, we have

\pro{Theorem 1.2}
Let $m=\deg I_{|P_{\xi}}$. There are constants $C_1$\ and $\eps_0$\
depending on $(r,\Del)$\ such that for any ample
divisor $H\sim a\Sigp+bP_{\xi}$\ with $b/a<\eps_0$\ and
any fiber type $\alp\ne 1_m$, we have
$$\mod\Rr(\alp)\leq (2r-1)d+C_1.
$$
\endpro

The proof of theorem 1.2 goes as follows: Let $\alpha=(\alp_1^{\oplus r_1},
\cdots,\alp_n^{\oplus r_n})$\ be any fiber type. Then each $E\in
\Rr(\alp)$\ admits a relative Hardar-Narasimhan filtration
$$0=E_0\sub E_1\sub\cdots\sub E_n=E\tag 1.4
$$
of which the quotient sheaves $F_i=E_i/E_{i-1}$\ are
torsion free with generic fiber types $(\alp_i^{\oplus
r_i})$\ respectively. Clearly, the deformation of $E$\ within $\Rr(\alp)$\
depends on deformation of individual $F_i$\ and the extension
$E_i\to E_{i+1}\to F_{i+1}$. The contribution of these data to the
number of moduli of $\Rr(\alpha)$\ can be estimated
by using Riemann-Roch. The details of the proof will be provided shortly.

In light of theorem 1.2, to prove theorem 1.1 we only need to study the
stratum $\Rr(1_m)$\ and
$$\Rr(1_m,D)=\{E\in\Rr(1_m)\mid \Hom(E,E( D))^0\ne0\}.
\tag 1.5
$$
In this section, we will first establish theorem 1.1 for
the stratum $\Rr(1_0,D)$\
and derive the remainder by induction on $m$.

Let $E$\ be any vector bundle of generic fiber type $1_0{=}(0^{\oplus r})$.
Let $x\in\Sig$\ be any point, let $P_x$\ be the fiber of $\pi$\ over
$x\in\Sig$\ and let $\beta_x(E)=(\beta_1^{\oplus r_1},
\cdots,\beta_n^{\oplus r_n})$\ be the fiber type of $E_{|P_x}$. In
case $\beta_x(E)\ne 1_0$, we call $P_x$\ a jumping line of $E$.
Let $P_x$\ be a jumping line of $E$. We then perform semistable reduction
on $E$\ along $P_x$\ by taking $F$\ to be the kernel of the (unique
surjective) homomorphism $E\to \OO_{P_x}(\beta_n)^{\oplus r_n}$.
For convenience, we will use $\Gam_x$\ to denote this operation and
denote $F=\Gam_x(E)$\ and $\ome_x(E)=\beta_n^{\oplus r_n}$. Clearly, $F$\
belongs to the exact sequence
$$0\lra F\lra E\mapright{\varphi} \OO_{P_x}(\beta_n)\uplu{r_n}\lra 0.
\tag 1.6
$$
An easy calculation based on Riemann-Roch yields

\pro{Lemma 1.3}
Let $F=\Gam_x(E)$\ with $\ome_x(E)=t^{\oplus s}$, then
$c_1(F)=c_1(E)-s[P_x]$\ and $c_2(F)=c_2(E)+s\cdott t$. In particular,
$\eta(r, c_2(F),c_1(F))=\eta(r, c_2(E),,c_1(E))+2rs\cdott t$.
\endpro

\proof See [Br, p166].\endpf

In case $F$\ still has a jumping line, say $P_{y}$\ of type
$(\cdots,\gamma_l^{\oplus s_l})$,
then we can further perform semistable reduction on $F$\ to get
$F_2=\Gam_y(F)$.
We can iterate this process as long as the resulting vector bundle $F_k$\
still admits jumping lines. In general, if $F_k$\ is derived by successively
performing this type of elementary transformations, namely, $F_0=E$\
and
$F_{i+1}=\Gam_{x_i}(F_i)$\ with $\ome_{x_i}(F_i)=t_i^{\oplus s_i}$\ for
$i=0,\cdots,k-1$, then we will write
$$F_k=\fretl(E), \quad
\Lambda=<x_1,\cdots, x_k>
$$
and define $\ome_{\Lambda}(E)=<t_1^{\oplus s_1},\cdots,t_k^{\oplus s_k}>$.
We call $k$\ the length of $\Lambda$.

\pro{Lemma 1.4}
For any vector bundle $E$\ of generic fiber type $1_0$, there is a finite
length $\Lambda=<x_1,\cdots,x_k>$\ such that $\fretl(E)$\ has no jumping
lines.
\endpro

\Proof
By lemma 1.3, the second Chern class
of $\Gam_x(E)$\ is strictly less that $c_2(E)$\ because $\beta_n<0$\
when $\beta_x(E)\ne 1_0$. Thus lemma 1.4 follows if we can show that
any vector bundle of generic fiber type $1_0$\
has non-negative second Chern class.
Indeed, let $E$\ be any vector bundle of generic fiber type $1_0$.
We choose a divisor $D$\ supported on fibers of $\pi$\ such that
$\OO_D$\ is a subsheaf of $E$\ with $E/\OO(D)$\ torsion free.
Since $E/\OO_D$\ has generic fiber type $1_0$,
we can assume $c_2(E/\OO(D))\geq 0$\ by the induction
hypothesis on the rank of $E$. Hence,
$$c_2(E)=c_2\Bl E/\OO(D)\Br +D\cdot(c_1(E)-D)=c_2\Bl E/\OO(D)\Br\geq 0.
$$
This completes the proof of lemma 1.4.
\endpf

Let $E$\ be a vector bundle of generic fiber type $1_0$\ and let
$\Lambda=<x_1,\cdots,x_k>$\ be such that $F=\Gam_{\Lambda}(E)$\ has
no jumping lines. Then $F$\ is a pull-back vector bundle $\pi\sta V$\
whose dual belongs to the exact sequence
$$0\lra E\dual\lra\pi\sta V\dual\lra J\lra 0,\tag 1.7
$$
where $J$\ is a torsion sheaf supported on the union of fibers $P_{x_i}$.
Usually, the sheaf $J$\ near some fiber $P_{x_i}$\ can be
very complicated.
The case that is easy to understand and will be dealt extensively
in the subsequent discussion is when $J\cong\OO_{P_{x_i}}(1)$\ near
$P_{x_i}$. The following theorem says that when the number of
moduli of $\Rr(1_0)$\ is close to $\eta(r,d,I)$, then for general $E\in
\Rr(1_0)$\ with the exact sequence (1.7), $J\cong\OO_{P_{x_i}}(1)$\
near $P_{x_i}$\ for all $x_i\in\{x_1\cdots x_k\}$\
except for a bounded number of fibers.

\pro{Theorem 1.5}
For any constant $e$, there is a constant $C_2$\ such that
$$\mod\Rr\leq \eta(r,d,I)+C_2.\tag 1.8
$$
Further, for any constant $C$, there are integers $l$, $l_1$, $l_2$\ and
$N_1$\ of which the following holds: Assume $d\geq N_1$\ and that
$S$\ is a variety parameterizing
a subset of $\Rr(1_0)$\ satisfying
$\mod(S)\geq \eta(r,d,I)-C$.
Then there is a line bundle $L$\ on $\Sig$\ of degree $[(d-c)/r]+l_1$,
where $c=I\cdot\sigp$, so that for general $E\in S$, there are
\roster
\item $d-l$\ distinct points $x_1,\cdots,x_{d-l}\in\Sig$\ in general
position, a surjective homomorphism $\tau_1\mh\pi\sta L\uplu{r}\to
\oplus^{d-l}_{i=1}\OO_{P_{x_i}}(1)$\ and
\item a zero dimensional scheme (divisor) $z_0\sub \Sig$\ away from
$\{x_1,\cdots,x_{d-l}\}$\ with $\len(z_0)\leq
l_2$\ and a sheaf of $\OO_{\pi^{-1}(z_0)}$-modules $J$\ with a quotient
homomorphism $\tau_0\mh\pi\sta L\uplu{r}\to J$\ so that $E\dual$\ belongs
to the exact sequence
$$0\lra E\dual\lra \pi\sta(L\uplu{r})\mapright{\tau_0\oplus \tau_1}
J \oplus \Bigl(\oplu{d-l}{i=1}\OO_{P_{x_i}}(1)\Bigr)\lra 0.\tag 1.9
$$
\endroster
\endpro

This theorem holds for a very simple reason: To maximize the number of
moduli of the set of those $E$\ in (1.9), we need to maximize the
number of moduli of the set of homomorphisms
$\tau_0\oplus\tau_1$\ and the quotient sheaves in (1.9). This can only
be achieved by letting $J=\{0\}$\ and $x_i$\ general.
Hence, if $\mod\Rr(1_0)$\ is close to the expected
dimension $\eta(r,d,I)$, then the number of fibers in
$\supp(J)$\ can not be too large.

Now we sketch how this structure theorem of $\Rr(1_0)$\ leads to the
proof of theorem 1.1. We first prove the case $m=0$\ by contradiction.
Assume $\mod\Rr(1_0,D)\geq \eta(r,d,I)-C$. Then by theorem 1.5, general
element $E\in\Rr(1_0,D)$\ fits into exact sequence (1.9) with $\{x_1,\cdots,
x_{d-l}\}$\ and $\tau_1\mh\pi\sta L^{\oplus r}\to\oplus\OO_{P_{x_i}}(1)$\
general. Now let $F=\ker\{\pi\sta L^{\oplus r}\mapright{\tau_0}J\}$\
and let $f\mh E\to E(D)$\ be a non-trivial traceless homomorphism. Then
 $f$\ induces a non-trivial traceless homomorphism
$$\tilde f: F\lra F(D+\pi\upmo(z_0)),
$$
where $z_0$\ is a divisor of $\Sigma$\ as in theorem 1.5 2).
Because the position of
$x_1,\cdots,x_{d-l}$\ and the homomorphism $\tau_1$\ are general,
we will see by degeneration theoretic methods
that for the torsion free sheaf $F\pri=\ker\bigl(
\pi\sta L^{\oplus r}\to \oplus(\OO_{P_{x_i}}\oplus
k_{p_i})\bigr)$, $\Hom(F\pri,F\pri(D+\pi\upmo(z_0)))\ne0$,
where $p_i\in P_{x_i}$.
Because of the special choice of $F\pri$, the non-vanishing of the
previous group amounts to say that for any choice of $p_i\in \Del$,
there are sections of $H^0(\OO(D+\pi\upmo(z_0)))$\ that vanishes on
$[(d-l)/r]$\ of $p_1,\cdots,p_{d-l}$. On the other hand,
since $D$\ is fixed and $\pi\upmo(z_0)$\ is bounded, this is impossible
if $\{p_i\}$\ are generic and $d$\ is sufficiently large. This leads to
a contradiction which ensures that $\eta(r,d,I)-\mod\Rr(1_0,D)$\ can
be arbitrary large.

For the general case, we use induction on $m$\ (with $r\geq m$ fixed).
Assume the theorem holds for $m-1\geq 0$\ and assume
$\mod\Rr(1_m,D)\geq\eta(r,d,I)-C$. Then for general $E\in\Rr(1_m,D)$, we can
perform an elementary transformation on $E$\ along a section $\Sigp$\
to get a new vector bundle $\tilde E\in \frak R^{r,d\pri}_{e\pri,I\pri}
(1_{m-1},\tilde D)$. By carefully study this correspondence,
we will get the
desired estimate of $\mod\Rr(1_m,D)$\ from the known estimate of $\mod
\frak R^{r,d\pri}_{e\pri,I\pri}(1_{m-1},\tilde D)$,
thus establishing the theorem 1.1.

In the following, we will fill in the details of the above sketch.
We continue
to use the notation introduced before lemma 1.4. We begin
with the estimate of the number of moduli of vector bundles of
generic fiber type $1_0$. Let $E_0\in\Rr(1_0)$\
be any vector bundle of generic fiber type $1_0$\ and $\Lambda
=<x_1,\cdots,x_k>$\ be such that $F=\Gamlam(E_0)$\ has no jumping
line. Then $F$\ is a pull-back vector bundle $\pi\sta V$. Let
$\ome$\ be $\ome_{\Lambda}
(E_0)=<t_1^{\oplus s_1},
\cdots,t_k^{\oplus s_k}>$\ and let
$$\SS_{\Lambda,\ome}(F)=\{E\in \Rr(1_0)\mid \ome\llam(E)=\ome\
\text{and}\ \Gam\llam(E)=F\}.
$$
In the following, we will estimate the number of moduli of this set.
We first study the case where $\Lambda=<x>$\ and $\ome=<t^{\oplus s}>$.
Let $\beta_x(F)=(\cdots,\beta_l^{\oplus r_l})$. Because of the following
lemma, either $t<\beta_l$\ or $t=\beta_l$.

\pro{Lemma 1.6}
Suppose $E_{|P_x}$\ has fiber type
$(\cdots,\gamma_n^{\oplus s_n})$\ and that $\Gam_x(E)$\
has fiber type $(\cdots,\beta_l^{\oplus r_l})$\
at $x$, then either $\gamma_n<\beta_l$\ or $\gamma_n=\beta_l$\
and $r_l\leq s_n$.
\endpro

\Proof
Since $F=\Gam_x(E)$\ is the kernel of
$E\to \OO_{P_x}(\gamma_n)\uplu{s_n}$,
$F_{|P_x}$\ belongs to the exact sequence
$$0\lra \OO_{P_x}(\gamma_n)^{\oplus s_n}\lra F_{|P_x}\lra \oplu{n-1}{i=1}
\OO_{P_x}(\gamma_i)^{\oplus s_i}\lra 0.
$$
Then the lemma follows because $\gamma_n<\gamma_{n-1}< \cdots<\gamma_1$.
\endpf

Let $E\in\SS_{x,\ome}(F)$. By dualizing the sequence (1.6), we get
$$0\lra E\dual \lra F\dual\lra \OO_{P_x}(-t)\uplu{s}\lra 0.\tag 1.10
$$
Clearly, all possible $E\dual$\ that fit into (1.10) are parameterized
by a subset of
$\Xi$\ that is the total space of $\Hom(F\dual,\OO_{P_x}(-t)\uplu{s})$.
Now let $\Theta\sub \Xi$\ be the subset
consisting of $\gamma\mh F\dual\to
\OO_{P_x}(-t)\uplu{s}$\ such that $\ker(\gamma)^{\vee}
\in \SS_{x,\ome}(F)$. $\Theta$\ admits a left $GL(s,\CC)$\ action and a right
$\Aut
(F\dual)$\ action as follows:
Let $\varphi_1\in \Aut(F\dual)$\ and let $\varphi_2\in GL(s)=
\Aut(\OO_{P_x}(-t)\uplu{s})$, then
$$\varphi_2\cdot\gamma\cdot\varphi_1=\varphi_2\circ\gamma\circ\varphi_1
\in \Hom(F^{\vee},\OO_{P_x}(-t)\uplu{s}).
$$
Geometrically, $\varphi\cdot\gamma\cdot\varphi_1$\ corresponds to a locally
free sheaf $E\pri$\ defined by
$$0\lra E^{\prime\vee}\lra F\dual \mapright{\varphi_2\cdot\gamma\cdot\varphi_1}
\OO_{P_x}(-t)\uplu{s}\lra 0.
$$
Clearly, $E\pri$\ is isomorphic to $E=\ker(\gamma)\dual$.
Conversely, suppose $E_1$\ and $E_2$\ are two isomorphic locally free
sheaves associated to $\gamma_1,\,\gamma_2\in\Theta$.
Then isomorphism $\varphi\mh E_1\to E_2$\ induces
isomorphism between $\Gam_x(E_1)$\ and $\Gam_x(E_2)$. Hence, there is an
automorphism $\varphi_1\mh F\dual\to F\dual$\ fitting into the
(commutative) diagram
$$\CD
0 @>>> E_2\dual @>>> F\dual @>\gamma_1>> \OO_{P_x}(-t)\uplu{s} @>>> 0 \\
@.  @VV\varphi\dual V @VV\varphi_1 V @. @. \\
0 @>>> E_1\dual @>>> F\dual @>\gamma_2>> \OO_{P_x}(-t)\uplu{s} @>>> 0 \\
\endCD
$$
In particular, there is a
$\varphi_2\mh\OO_{P_x}(-t)\uplu{s}\to\OO_{P_x}(-t)\uplu{s}$\ such that
$\varphi_2\circ\gamma_1=\gamma_2\circ\varphi_1$.
Therefore those points in $\Theta$\ that give rise to
isomorphic sheaves form an $\Aut(F\dual)\times
GL(s)$\ orbit. Next, we will determine the size of the stabilizer in
$\Aut(F\dual)\times GL(s)$\ of any $\gamma\in\Theta$. Suppose
$\varphi_1\in\Aut(F\dual)$\ and $\varphi_2\in GL(s)$\ are such that
the right rectangle below is commutative,
$$\CD
0 @>>> E\dual @>>> F\dual @>\gamma>> \OO_{P_x}(-t)\uplu{s} @>>> 0 \\
@. @. @VV\varphi_1 V   @VV\varphi_2^{-1} V @. \\
0 @>>> E\dual @>>> F\dual @>\gamma>> \OO_{P_x}(-t)\uplu{s} @>>> 0 \\
\endCD
$$
Then it induces a $\varphi\in\Aut(E\dual)$. One sees
that such a map $\text{Stab}_{\gamma}\to \Aut(E\dual)$\ is injective. Thus,
if we let $\Theta_k\sub\Theta$\ be the set of $\gamma$'s such that
$\aut(\ker(\gamma))=k$\ (for a set $R$\ of sheaves, we define
$\aut(R)=\max_{E\in R}\{\dim\Hom(E,E)\}$), then
$$\dim\BBl\text{GL}(s)\backslash\Theta_k/\text{Aut}(F\dual)\BBr
\leq\dim\Hom(F\dual,\OO_{P_x}(-t)^{\oplus s})-
(s^2+\aut(F\dual))+k.\tag 1.11
$$
Finally, let $(\beta_1^{\oplus r_1},\cdots,\beta_n^{\oplus r_n})$\
be the fiber type of $F_{|P_z}$, then by lemma 1.6, $t\leq \beta_i$.
Because $\sum r_i\beta_i=0$, we have
$$\dim \Hom\Bl F\dual,\OO_{P_x}(-t)\Br=\sum^n_{i=1} r_i\dim H^0(\OO_{P_x}
(\beta_i-t)=r(-t+1).\tag 1.12
$$

Returning to the general case $\Lambda=<x_1,\cdots,x_k>$\ and $\ome=
<t_1^{\oplus s_1},\cdots,t_n^{\oplus
s_n}>$, we will prove:

\pro{Lemma 1.7}
With the notation as above and let $E\in\sslo(F)$, then
$$\align
(\mod&-\aut)\Bl\sslo(F)\Br\leq\\
&\leq\eta(E)-\BBl r\sum^n_{i=1}s_i(-t_i-1)+ \sum^n_{i=1}s_i^2\BBr
-\aut(F)-(r^2-1)(g-1).\tag 1.13
\endalign
$$
\endpro

\proof
We only need to prove the inequality
$$(\mod-\aut)\Bl\sslo(F)\Br\leq \sum^{n}_{i=1}\Bl rs_i(-t_i+1)-s_i^2\Br
-\aut(F)\tag 1.14
$$
because then (1.13) follows from
$c_2(E)=-\sum_{i=1}^n s_it_i$\ and $\eta(E)=-2r\sum^n_{i=1}s_it_i+
(r^2-1)(g-1)$. We prove (1.14) by induction on $n$. When $n=1$, (1.14)
follows from (1.11) and (1.12) because
$\aut({\Cal S}_{\Lambda,\ome}(F))=\sup\{k|\Theta_k\ne\emptyset\}$.
Now assume (1.14) is true for $n-1$. We divide $\SS_{x_1,\ome_1}(F)$,
$\ome_1=(t_1\uplu{s_1})$, into subsets $W_k$\ such that $F\pri\in W_k$\
if $\aut(F\pri)=k$. Let $\Lambda_2=<x_2,\cdots,x_k>$\ and
$\ome_2=<t_2\uplu{s_2},\cdots,t_n \uplu{s_n}>$.
Then by induction hypothesis, for $F\pri\in W_k$,
$$(\mod-\aut)\Bl \SS_{\Lambda_2,\ome_2}
(F\pri)\Br\leq\sum_{i=2}^n \Bl rs_i(-t_i+1)-s_i^2\Br-k
$$
and therefore,
$$\align
(\mod-\aut)(\sslo(F))&\leq\sup_k\biggl\{
\sum_{i=2}^n\Bl rs_i(-t_i+1)-s_i^2\Br-k+\mod(W_k)\biggr\} \\
&\leq \sum_{i=2}^n\Bl rs_i(-t_i+1)-s_i^2\Br + \BBl rs_1(-t_1+1)-s_1^2-
\aut(F)\BBr\\
&=\sum_{i=1}^n\Bl rs_i(-t_i+1)-s_i^2\Br -\aut(F).\tag{\qed}
\endalign
$$

Now we are ready to prove our structure theorem for subsets of $\Rr(0_r)$.

{\it Proof of theorem 1.5}.
(1.8) follows directly from Riemann-Roch and the fact that there is a constant
$C_2$\ depending on $(\Del,H,e)$\ such that for any $E\in \Rr$,
$\aut (E)\leq C_2$.
We now prove the second part of the theorem.
Let $S\sub\Rr(1_0)$\ be any (irreducible) algebraic set and let
$E\in S$\ be a generic element. By lemma 1.4, after performing
a sequence of semistable reduction at $y_1,\cdots,y_n$, we get a vector
bundle with no jumping line, say $\pi\sta F$\ with $F$\ a vector bundle
over $\Sig$. Clearly, $n=n(E)$\ depends on $E$. We let $S_0\sub S$\ be
the open set of $E\pri\in S$\ with $n(E\pri)=E$\ and let $n_0$\
be the integer so that when $E$\ varies in $S_0$,
the number of moduli of the (unordered) set
$y_1,\cdots,y_n$\ is $n_0$. In other words, $n_0$\ of $(y_1,\cdots y_n)$\
are in generic position. We know that the number of moduli of rank $r$\
vector bundles on $\Sig$\ is $r^2(g-1)+1$. Also,
since $E$\ is $e$-stable, $\aut(E)$\ is bounded by a constant
$C_3\pri$\ independent of $d$\ and $I$\ (see lemma 1.10). Combining these
with (1.13), we get
$$\mod(S)\leq\eta(E)-r\sum^n_{i=1}s_i(-t_i-1)-\sum^n_{i=1}s_i^2+n_0
-\aut(F)+g+C_3\pri.\tag 1.15
$$
Since we have assumed $\mod(S)\geq \eta(r,d,I)-C$, for $C_3=C+C_3\pri+g$,
we get
$$C_3\geq r\sum^n_{i=1}s_i(-t_i-1)+\BBl\sum^n_{i=1}s_i^2-n_0\BBr
+\aut(F).\tag 1.16
$$
Because $t_i<0$, all terms in (1.16) are non-negative. This immediately
gives us $n-n_0\leq C_3$. Next,
we define the multiplicity $m(y_i)$\ of $y_i$\ to be the number of appearance
of the point $y_i$\ in $(y_1,\cdots,y_n)$. Then by (1.16),
$${1\over 2}\#\{y_i\mid m(y_i)\geq 2\}\leq \sum^n_{i=1}s_i^2-n_0\leq
C_3.
$$
So the total multiplicity of multiple points is bounded.
Without loss of generality, we can assume $y_1,\cdots, y_{n_0}$\ are in general
position for general $E\in S$. For convenience, we call
$y_i\in (y_1,\cdots,y_n)$\ a simple point if $m(y_i)=1$\ and $\ome_{y_i}(E)=
t_i^{\oplus s_i}$\ is $(-1)^{\oplus 1}$. We claim that then
$$\sum_{y_i\ne \text{simple}}(-s_it_i)\leq 2\sum^n_{i=1}s_i(-t_i-1)+
(\sum_{i=1}^n s_i^2-n_0)+
\#\{y_i \mid m(y_i)\geq 2\}\leq ({2\over r}+2)C_3.\tag 1.17
$$
Indeed, when $t_i\leq -2$, then the term $s_i(-t_i)$\ is bounded
from above by term $2s_i(-t_i-1)$\ in the middle of (1.17) and when $t_i=-1$\
and
$s_i\geq 2$, then we have $-s_it_i\leq s_i^2-1$.
The only remaining situation is when $m(y_i)\geq 2$, $t_i=1$\
and $s_i=1$. But in this case, $(-t_i)s_i=1$\ can be absorbed by term
$\#\{y_i\mid m(y_i)\geq 2\}$. Hence, (1.17) holds.
Finally, since $d=\sum^n_{i=1}(-t_i)s_i$,
$$\#\{y_i\mid y_i\ \text{simple}\}=d+\sum_{y_i\ne \text{simple}}s_it_i
\geq d-4C_3.
$$
Therefore, combined with $n-n_0\leq C_3$, we get
$$d\geq n\geq n_1=\#\{y_i\mid y_i\ \text{simple}, 1\leq i\leq n_0\}\geq
d-5C_3.  \tag 1.18
$$
Now we let $l=[5C_3]+1$.
Without loss of generality, we can assume $\{y_1,\cdots,y_{d-l}\}$\
are simple points in $\{y_1,\cdots,y_{n_0}\}$.
Then the sheaf $E$\ must belong to the exact sequence
$$0\lra \pi\sta F\lra E\lra \Bigl(\oplu{d-l}{i=1}\OO_{P_{y_i}}(-1)\Bigr)
\oplus J\pri\lra 0.\tag 1.19
$$
To prove the proposition, we need to have an estimate on $F$\ and $J\pri$.
By definition, $J\pri$\ admits a filtration
$$0=J_{d-l}\sub J_{d-l+1}\sub\cdots\sub J_n=J\pri
$$
such that $J_{i+1}/J_i\cong \OO_{P_{y_i}}(t_i)\uplu{s_i}$. Thus there is
a zero scheme $z\pri\sub \Sig$\ supported on $\{y_{d-l+1},\cdots,y_n\}$\
of length $\len(z\pri)\leq n-(d-l)\leq 5C_3$\ (because of (1.18))
such that $J\pri$\ is an $\OO_{\pi^{-1}(z\pri)}$-modules and further
$$0\leq c_1(E)\cdot\Sigp -(\deg F+d-l)= c_1(J\pri)\cdot \sigp=
\sum_{i=d-l+1}^n s_i \leq \sum_{i=d-l+1}^n(-t_i)s_i
\leq 5C_3.\tag 1.20
$$
Here, the last inequality holds because of (1.17) and $n-n_0\leq C_3$.
Also, since $\aut (F)\leq C_3$\ (from (1.16)),
there is a constant $C_4$\ such that $F$\ is $C_4$-stable.

It remains to show that we can find an integer $l_1$\ (independent of $d$)
and find a single line bundle $L$\ of degree $[(d-c)/r]+l_1$\ ($c=I\cdot\sigp$)
so that for any $E\in S_0$, $E$\ belongs to the exact sequence
$$0\lra E\dual \lra \pi\sta(L\uplu{r})\lra \Bigl(\oplu{d-l}{i=1}
\OO_{P_{y_i}}(1)\Bigr) \oplus J\lra0\tag 1.21
$$
specified in theorem 1.5.
First, there is a constant $l_1$\ and a line bundle $L$\ of degree $
[(d-c)/r]+l_1$\ such that for any $C_4$-stable rank $r$\ vector bundle
$F$\ on $\Sig$\ satisfying (1.20), $L\otimes F$\ is generated by
$H^0(L\otimes F)$. Now for any $E\in S_0$\ with the data given by
(1.19), we choose $\pi\sta F\dual\to \pi\sta L\uplu{r}$\
so that the support of
$\pi\sta(L\uplu{r})/\pi\sta F\dual$\ is
disjoint from $\cup_{i=1}^{d-l} P_{y_i}$. Then
by dualizing (1.19) and coupled with
$\pi\sta F\dual\to \pi\sta(L\uplu{r})$, we get
$$0\lra E\dual\lra \pi\sta(L\uplu{r})\lra J\oplus \Bigl(\oplu{d-l}{i=1}
\OO_{P_{y_i}}(1)\Bigr) \lra 0.
$$
Finally, it is easy to see that there is an integer $l_2$\ depending
only on $C_3$\ and $l_1$\ such that for some subscheme $z\sub\Sigma$\
of length $\len(z)\leq l_2$, $J$\ is a sheaf of $\OO_{\pi^{-1}(z)}$-modules.
This completes the proof of the theorem.
\endpf

Now we prove theorem 1.2.

{\it Proof of theorem 1.2}.
We begin with a general vector bundle $E\in \Rr(\alp)$,
$\alp=(\alp_1\uplu{r_1},
\cdots,\alp_n\uplu{r_n})\ne 1_m$. Let
$$0=E_0\sub E_1\sub\cdots\sub E_n=E\tag 1.22
$$
be the relative Hardar-Narasimhan filtration
such that $F_i=E_i/E_{i-1}$\ are torsion free of
generic fiber types $(\alp_i\uplu{r_i})$\ respectively. We
call this the relative filtration of $E$. ((1.12) can be derived by using
the usual Harder-Narasimhan filtration of $E$\ with respect to the
divisor $kP_{\xi}+\Sigp$\ with $k\gg0$.)

We fix $F_i=E_i/E_{i-1}$\ and let $W(\{F_i\}^n_1)$\
be the set of all vector bundles $V$\ such that whose relative
filtrations $0\sub V_1\sub\cdots\sub V_n=V$\ satisfy $V_i/V_{i-1}
\cong F_i$. Our first step is to estimate the number of local
moduli $\modl\Bl W(\{F_i\}^n_1)\ \text{at}\ [E]\Br$.
Let $A_i={1\over r_i}c_i(F_i)$\
and $d_i=c_2(F_i)-({r\atop 2})A_i^2$. Note that by proof of lemma 1.4,
$d_i\geq0$. Now an easy calculation shows that
$$d=c_2(E)=\half\sum_{i=1}^n r_i(I-A_i)\cdott A_i+\sum_{i=1}^n d_i
=\half I^2-\half\sum_{i=1}^n r_i A_i^2+\sum_{i=1}^n d_i.\tag 1.23
$$
{}From the exact sequence
$$0\lra E_{n-1}\lra E\lra F_n\lra 0
$$
and the argument similar to (1.11), we have
$$\align
\modl\Bl W(\{F_i\}^n_1)\, &\text{at}\,[E]\Br\leq
\modl\Bl W(\{F_i\}^{n-1}_1)\,\text{at}\,[E_{n-1}]\Br+
\dim \Ext^1(F_n,E_{n-1})\\
&-\aut(F_n)-\aut(E_{n-1})+\aut(E)-\dim\Hom(F_n,E_{n-1}).
\tag 1.24
\endalign
$$
Further, because $E_{n-1}\dual\otimes F_n$\ has generic fiber type
$((\alp_n-\alp_1)^{\oplus r_1},\cdots,(\alp_n-\alp_{n-1}^{\oplus r_{n-1}}))
$\ and $\alp_{i+1}<\alp_i$, $\Ext^2(F_n,
E_{n-1})=0$\ by Serre duality. Hence
$$\dim\Hom(F_n,E_{n-1})-\dim\Ext^1(F_n,E_{n-1})=\xx(F_n, E_{n-1}),
\tag 1.25
$$
where the right hand side of (1.25) is the abbreviation of
$\xx(\Ext^{\cdot}(F_n,E_{n-1}))$. Finally, by using the filtration (1.22),
we have
$$\align
\modl\Bl & W(\{F_i\}^n_1)\,\text{at}\,[E]\Br-\aut(E)\\
&\leq \modl\Bl W(\{F_i\}^{n-1}_1)\,\text{at}\,[E_{n-1}]\Br-
\sum_{i=1}^{n-1}\xx(F_n,E_{n-1})-\aut(E_{n-1})-\aut(F_n)\\
&\leq \sum_{i>j}\xx(F_i,F_j)-\sum_{i=1}^n\aut(F_i).
\tag 1.26
\endalign
$$
The last inequality is derived by iterating the first part of (1.26).
Therefore,
$$\align  \mod\Rr(\alp)\leq&\sup\biggl\{
 \sum_{i>j}\xx(F_i,F_j)
+\sum_{i=1}^n\bigl(\mod(F_i)-\aut(F_i)\bigr)\biggr\}\\
&+\max\{\aut(E)\mid E\in\Rr(\alp)\},
\tag 1.27
\endalign
$$
where the sup is taken over all possible relative
filtrations (1.22) of $E$'s in $\Rr(\alpha)$. We now
calculate the right hand side of (1.27) by Riemann-Roch. First,
$$\xx(F_i,F_j)=
r_ir_j\BBl\half(A_j-A_i)^2-\half(A_j-A_i)\cdot\! K_{\Del}+(1-g)\BBr
-r_id_j-r_jd_i.
$$
For simplicity, in the following we will group all terms that are bounded
independently of $r_i,\,d_i$, $A_i$\ and $\alp_m\ne1_m$\
into $O(1)$.
We have
$$
\sum_{i>j}\xx(F_i,F_j)=-\sum_{i>j}r_ir_j\Bigl(\half(A_j-A_i)^2-\half(A_j-A_i)\cdot\!K_{\Del}
+{d_i\over r_i}+{d_j\over r_j}\Bigr) +O(1).\tag 1.28
$$
Further, one calculates
$$\align
\eta(F_i)&=2r_i\Bl d_i+\Bigl({r_i\atop 2}\Bigr)A_i^2\Br -
(r_i-1)r_i^2A_i^2-(r_i^2-1)(1-g) \\ & =2r_id_i-(r_i^2-1)(1-g).
\tag 1.29
\endalign
$$
Thus by combining (1.23), (1.27)-(1.29) and the fact that $\mod(F_i)-
\aut(F_i)\leq \xx(F_i,F_i)$, we obtain
$$\align
(\mod(E)&-\aut(E))\Bl\Rr(\alpha)\Br-(2r-1)d\\
&\leq-\sum_{i>j}r_ir_j\Bigl(\half(A_j-A_i)^2-\half(A_j-A_i)\cdott K_{\Del}
+{d_i\over r_i}+{d_j\over r_j}\Bigr)\\
&\quad\quad+\sum_{i=1}^n 2r_id_i -(2r-1)\Bigl(\sum_{i=1}^nd_i-{1\over 2}
\sum_{i=1}^nr_iA_i^2\Bigr)+O(1).\tag 1.30
\endalign
$$
To analyze (1.30), we first note that
$$\sum_{i>j}r_ir_j(A_j-A_i)^2=r\sum_{i=1}^n r_iA_i^2-I^2;
\quad \sum_{i>j}r_ir_j\Bigl({d_i\over r_i}+{d_j\over r_j}\Bigr)=
\sum_{i=1}^n (r-r_i)d_i.
$$
Now if we let $A_i\nue \alp_i\sigm+c_i\Px$\ and let
$\delta=\sigp\cdott\sigp$, then the right hand side of (1.30) is equal to
$$\align
\sum_{i=1}^n\Bigl(\half(r-1) & r_iA_i^2-
(r-r_i-1)d_i\Bigr)+\half\sum_{i>j}r_ir_j(A_j-A_i)\cdot\! K_{\Del}+O(1)\\
&=\sum_{i=1}^n\half(r-1)r_i(-\delta\alp_i^2+2\alp_i c_i)
-\sum_{i=1}^n (r-r_i-1)d_i\\
&\quad+ \half\sum_{i>j}r_ir_j(\alp_j-\alp_i)\sigm\cdot\! K_{\Del}+
\half\sum_{i>j}r_ir_j(c_j-c_i)\Px\cdot\! K_{\Del}+O(1)
\endalign
$$
which is bounded from above by (note $d_i\geq 0$)
$$\align
-{1\over 4}\delta& \sum_{i=1}^n\alp_i^2+(r-1)\sum_{i=1}^nr_i\alp_i c_i-
\sum_{i>j}r_ir_j(c_j-c_i)+O(1)\\
=&-{1\over 4}\delta\sum_{i=1}^n\alp_i^2+
\sum_{k=1}^n\biggl(r_kc_k\Bigl((r-1)\alp_k+\sum_{i=1}^{k-1}r_i-
\sum_{i=k+1}^nr_i \Bigr)\biggr) +O(1).
\tag 1.31\endalign
$$
Let $p_k=(r-1)\alp_k+\sum_{i=1}^{k-1}r_i-\sum_{i=k+1}^nr_i$. Then
when $\alp\ne 1_m$, or equivalently when $n\geq 2$\ or $n=2$\
and $\alp_1-\alp_2\geq 2$, we will always have
$$p_k-p_{k+1}=(r-1)(\alp_k-\alp_{k+1})-(r_k+r_{k+1})\geq0,\quad k\leq n-1.
$$
We rewrite
$$\sum_{k=1}^n\biggl(r_kc_k\Bigl((r-1)\alp_k+\sum_{i=1}^{k-1}r_i-
\sum_{i=k+1}^nr_i \Bigr)\biggr)
= \sum_{k=1}^{n-1}\biggl((p_k-p_{k+1})\Bl\sum_{i=1}^kr_ic_i\Br\biggr)
+p_n\sum_{i=1}^nr_ic_i.
$$
Finally, we shall make use
the fact that $E$\ is $e$-stable. If $H\nue a\sigp+b\Px$, then for
any $k$,
$$\deg(E_k)=a\sum_{i=1}^k r_ic_i+b\sum_{i=1}^k r_i\alp_i\leq
{\rk (E_k)\over r}H\cdott I+e\sqrt{H^2}.
$$
Therefore, for $k\leq n-1$,
$$ \sum_{i=1}^k r_ic_i \leq 1+{b\over a} +e\sqrt{\delta+2{b\over a}}-
{b\over a} \sum_{i=1}^kr_i\alp_i.\eqno
$$
Thus we get
$$\align
\sum_{k=1}^n \Bigl((p_k-&p_{k+1})(\sum_{i=1}^k r_ic_i)\Bigr)\\
&\leq \sum_{k=1}^{n-1} (p_k-p_{k+1})
\Bigl({a+b+e\sqrt{a^2\delta+2ab}\over a}-{b\over a}\sum_{i=1}^k r_i\alp_i
\Bigr) +p_n m\\
&\leq {b\over a}\cdot r^2(\sum_{i=1}^n|\alp_i|)^2+
4r^2(2+e\delta)(1+{b\over a}))(\sum_{i=1}^n|\alp_i|)+O(1).
\tag 1.32\endalign
$$

Here we have used the fact that $p_k-p_{k+1}\leq \sum_{i=1}^n|\alp_i|+r$\
and $p_n\leq 0$\ because $\alp\ne1_m$\ and $\sum_{i=1}^nr_ic_i=m\geq 0$.
Now if we assume
$${b\over a}r^2<{1\over 16}\delta,
$$
then everything in (1.32) can
be absorbed by the quadratic  term $-{1\over 4}\delta\sum_{i=1}^n\alp_i^2$\
(in (1.31)) with the help of some constant $C_1$.
Thus combined with (1.31), we have proved
$$\mod\Rr(\alp)\leq (2r-1)d+C_1+
\max\{\aut(E)\mid E\in\Rr(\alp)\}.
$$
The theorem 1.2 will be proved if we can bound $\Hom(E,E)$\
for $E\in\Rr$.
Since $E$\ is $e$-stable, $E\dual\otimes E$\ must be $2|e|+1$-stable.
(This can be proved by using the fact that the Harder-Narasimhan filtration
of $E$\ will induce the Hardar-Narasimhan filtration of $E\dual\otimes E$.)
Thus $\aut(E)$\ is bounded independently of $d$\ by the following lemma.

\pro{Lemma 1.8}
For constants $e_1$, $e_2$\ and integer $r$, there is a constant $C\pri$\ such
that
whenever $V$\ is a rank $r$\ $e_1$-stable vector bundle on $\Del$\ such that
$|\deg(V)|\leq e_2$, then we have
$\dim H^0(V)\leq C\pri$.
\endpro

\Proof
We prove the lemma by induction on $r$. The case $r=1$\ is obvious.
Assume the lemma is true for vector bundles of rank $\leq r-1$\ and assume $V$\
has $H^0(V)\ne\{0\}$. Then there is a line bundle $L$,
$\deg L\geq 0$\ such that $V$\ belongs to the exact sequence
$$0\lra L\lra V\lra V/L\lra 0
$$
with $V/L$\ torsion free. Since $V$\ is $e_1$-stable and $|\deg(V)|\leq
e_2$, there are constants $e\pri_1$\ and $e\pri_2$\ such that $|\deg L|$,
$|\deg V/L|\leq e\pri_2$\ and $V/L$\ is $e\pri_1$-stable. Thus
by induction hypothesis, there is a constant $C\pri$\ such that
$h^0(L)\leq C\pri$\ and $h^0(V/L)\leq C\pri$.
The lemma then follows.
\endpf


We now prove theorem 1.1 by induction on $m$.
We first establish the case $m=0$. Let $e$\ and $C$\ be any constants,
$r\geq 2$\ be an integer and $D\sub \Del$\ be any divisor. We assume $H$\
is an ample divisor satisfying the condition of theorem 1.2.
To prove the theorem, we need to show that there is a constant $N$\ depending
only on $(X,H,r,I,e,D)$\ so that if for some $d$\ we have
$$\mod\Wr(D)\geq \eta(r,d,I)-C,\tag 1.33
$$
then $d\leq N$. Now assume (1.33) does hold.
Thanks to theorem 1.2, there is an $N_1\geq0$\ such that
if $d\geq N_1$, then the set
$\Wr(D,1_0)$\ satisfies
$$\mod\Wr(D,1_0)=\mod\Wr(D)\geq \eta(r,d,I)-C.\tag 1.34
$$
Of course, $\Wr(D,1_0)$\ is a constructible set. Let $S$\ be
an irreducible variety parameterizing a subset of
$\Wr(D,1_0)$\ such that $\mod S\geq\eta(r,d,I)-C$.
By theorem 1.5, there are constants $l$, $l_1$, $l_2$\ (independent of
$d$) and line bundle $L$\ of degree $[(d-c)/r]+l_1$\
such that associated to a general $E\in S$, there
are $x_1,\cdots,x_{d-l}\in\Sig$\ in general position and a quotient sheaf
$J_E$\ of $\pi\sta(L\uplu{r})$\ such that $E$\ belongs to the exact sequence
$$0\lra E\dual\mapright{i} \pi\sta(L\uplu{r})\mapright{\tau_0\oplus \tau_1}
J_E\oplus \Bigl(\oplu{d-l}{i=1}\OO_{P_{x_i}}
(1)\Bigr)\lra0.\tag 1.35
$$
Clearly, $E$\ is determined by the surjective homomorphisms
$$\pi\sta\Bl L\uplu{r}\Br\mapright{\tau_0}J_E\ \text{and}\
\pi\sta\Bl L\uplu{r}\Br\mapright{\tau_1}\oplus\OO_{P_{x_i}}(1).
$$
Hence the combined number of moduli of the sets of
these quotient sheaves
that come from $E\in S$\ is no less than $\eta(r,d, I)-C$. Let
$$\Xi_0=\bigl\{\tau_0:\pi\sta\Bl L\uplu{r}\Br\to J_E\mid E\in S\bigr\}
$$
$$\Xi_1=\bigl\{\tau_1:\pi\sta\Bl L\uplu{r}\Br\to \oplu{d-l}{}\OO_{P_{x_i}}
(1)\mid E\in S\bigr\}.
$$
Because of the following lemma, the information contained in
$\Xi_0$\ is minimal.

\pro{Lemma 1.9}
There is a constant $C_5$\ independent of $d$\ such that
$\mod\Bl\Xi_0\Br\leq C_5$.
\endpro

\proof
We first calculate the Hilbert polynomials of the sheaves $J_E$. Let $J_E(n)=
J_E\otimes H^{\otimes n}$. Then
$$
\xx(J_E(n))=\xx(\pi\sta L\uplu{r}(n))-\xx(E(n))-(d-l)\xx(\OO_{P_{\xi}}(1)
\otimes H^{\otimes n})
=a_1(d)n+a_0(d),
$$
where $a_1(d)=(r[(d-c)/r]-d+rl_1+l)\cdott (H\cdott P_{\xi})-
I\cdott H$\ and $a_0(d)=(r[(d-c)/r]-d)+rl_1-\half I^2+\half
I\cdott K+2l$.
Since for integers $d$,
$r[(d-c)/r]-d$\ can only attain integer values between
$-c-r$\ and $-c$, the function $a_1(d)$\ (resp. $a_0(d)$)
can only attains $r$\ values. Hence, $\{\xx\Bl
J_E(\cdot)\Br \mid E\in S\}$\ is a finite set (independent of $d$)
and by [Gr, p12], the set $\Xi_0$\ is bounded. Thus, there is a
constant $C_5$\ such that $\mod\Xi_0\leq C_5$.
\endpf

Since $\mod S\leq \mod\Xi_0+\mod\Xi_1$, we have
$$\mod\Xi_1\geq\eta(r,d,I)-(C+C_5).
$$
Let $\tau_1\in\Xi_1$\ and $F=\ker\{\tau_1\}$.
In the following, we seek to relate the the
non-vanishing of $\Hom(E,E(D))^0$\ to the
non-vanishing of $\Hom(F,F(D\pri))^0$\ for some divisor $D\pri$.
First of all, by 2) of theorem 1.5, there is a divisor $z\in
\Sigma$\ (of degree $\leq l_2$) such that the composition
$$F(-\pi\upmo(z))\hookrightarrow\pi\sta L(-\pi^{-1}(z))
\uplu{r}\lra \pi\sta L\uplu{r}\mapright{\tau_0} J_E \tag 1.36
$$
is trivial. Because of (1.35),
$F(-\pi^{-1}(z))$\ is a subsheaf of $E^{\vee}$. Therefore, any
non-trivial traceless homomorphism $\varphi\mh E\to E(D)$\
will provide us a non-trivial traceless homomorphism
$$F(-\pi^{-1}(z))\lra E^{\vee}\mapright{\varphi} E^{\vee}(D)\lra F(D)
$$
Further, let $\bar z$\ be a fixed divisor on $\Sigma$\ of degree $l_2+2g$.
Since
$h^0(\Sigma,\OO_{\Sigma}(\bar z-z))\ne0$, $\Hom\Bl F,F
(D\plu\pi^{-1}(z)\Br\ne0$\ implies
$\Hom\Bl F,F(D\plu\pi^{-1}(\bar z))\Br\ne0$. Thus we have proved:

\pro{Lemma 1.10}
With the notation as before, then there is a divisor $z\subset \Sigma$\
independent of $d$\ and $D$\ such that for any sheaf $F=\ker\{\tau_1\}$,
where $\tau_1\in\Xi_1$, and for $D_1=D+\pi^{-1}(z)$, we have
$\Hom\Bl F,F(D_1)\Br\ne\{0\}$.\qqed
\endpro

Our next step is to investigate the set $\Xi_1$\ by utilizing this
non-vanishing property. We first fix $d-l$\ general
points $x_1,\cdots,x_{d-l}\in\Sigma$\ and let $U$\ be the set of all quotient
homomorphisms
$$\sigma: \pi\sta\Bl L\uplu{r}\Br\lra \oplu{d-l}{i=1}\OO_{P_{x_i}}(1).
\tag 1.37
$$
$U$\ is (canonically) parameterized by an open subset of the
product of $d-l$\ copies of projective space
$\PP^{2r-1}$\ after fixing basis of each $H^0(\OO_{P_{x_i}}(1))$.
In the following, for any $u\in \Pi^{d-l}\PP^{2r-1}$\ we denote by
$\sigma_u$\ the associated homomorphism $\sigma_u\mh \pi\sta(L\uplu{r})
\to \oplus^{d-l}_{i=1}\OO_{P_{x_i}}(1)$. Let
$$\Xix=\{u\in \Pi^{d-l}\PP^{2r-1}\mid \sigma_u\in\Xi_1\}.
$$
Since the points of $({\frak x})=(x_1,\cdots,x_{d-l})$\ are general,
$$\dim \Xix\geq \mod \Xi_1-(d-l)\geq \eta(r,d,I)-(C+C_5)-(d-l)\geq (2r-1)d-
C_6 \tag 1.38
$$
for some integer $C_6$. Now let $l_3=l+C_6+1$. Possibly after rearranging
the order of $(x_1,\cdots,x_{d-l})$,
we can further assume that the restriction to $\Xix$\ of the projection
from $\Pi^{d-l}\PP^{2r-1}$\ to the first
$d-l_3$\ factors,
$$\Xix\sub \Pi^{d-l}\PP^{2r-1}\lra\Pi^{d-l_3}\PP^{2r-1},
$$
is dominant. That is, for general $v\in \Pi^{d-l_3}\PP^{2r-1}
$\ with $\sigma_v\pri\mh\pi\sta L\uplu{r}\to \oplus^{d-l_3}_{i=1}\OO_{P_{x_i}}
(1)$\ the associated homomorphism, there is at least
one $\xi\mh \pi\sta L\uplu{r}\to \oplus^{d-l}_{i=d-l_3+1}\OO_{P_{x_i}}(1)$\
such that $\sigma\pri_v\oplus\xi$\ considered as a quotient sheaf
belongs to $\Xix$.
Thus if we let $V=\ker\{\sigma\pri_v\oplus\xi\}$\ and let
$V_v=\ker\{\sigma_v\pri\}$,
$V$\ and $V_v$\ fit into the following exact sequence
$$0\lra V\lra V_v\lra\oplu{d-l}{i=d-l_3+1}\OO_{P_{x_i}}(1)
\lra 0.
$$
Put $A=\cup_{i=d-l_3+1}^{d-l} P_{x_i}$\ be a divisor in $\Del$.
Following the argument in lemma 1.10, the non-trivial homomorphism $\phi$\
in lemma 1.10 induces a non-trivial homomorphism
$$\phi\pri\mh V_v\lra V_v(D_1+A).\eqno
$$
Therefore for general $v\in \Pi_{i=1}^{d-l_3} \PP^{2r-1}$,
$\Hom_{\Del}(V_v,V_v(D_1\plu A))^0\ne\{0\}$.
Finally, as in lemma 1.9, for any fixed divisor $A_0\sub \Del$\ consists
of $l_3+2g$\ fibers of $\Del$, we must have $\Hom\Bl V_v,V_v
(D_1+A_0)\Br^0\ne0$\ as well. Therefore, theorem 1.1 (when $m=0$) follows from

\pro{Proposition 1.11}
For any divisor $D\sub \Del$\ and any integer $l_0$, there is a constant
$N$\ of which the following holds:
Assume $d\geq N$, that $L$\ is a line bundle on $\Sig$\
of $\deg L=[d/r]+l_0$\ and that $x_1,\cdots,x_d$\ are general points in $\Sig$,
then for general $v\in \Pi^d\PP^{2r-1}$, the sheaf $E_v=\ker
\{\sigma_v\}$, where $\sigma_v$\ is the associated homomorphism
$\pi\sta L^{\oplus r}\to \oplus^d_{i=1} \OO_{P_{x_i}}(1)$, satisfies
$\Hom(E_v,E_v(D))^0=\{0\}$.
\endpro

\Proof
We prove it by contradiction. The trick is to first prove the
vanishing of this homomorphism group for a special quotient sheaf
and then apply the semicontinuity theorem to derive the general
case.
Let $p_i\in P_{x_i}$\ be general closed point
and let $U$\ be a small disk containing 0.
There is a torsion free sheaf $J_i$\ on
$P_{x_i}\timess U$\ flat over $U$\ such that $J_{i|P_{x_i}\timess\{0\}}\cong
\OO_{P_{x_i}}\oplus \CC_{p_i}$\
($\CC_{p_i}$\ is the skyscraper sheaf supported on $p_i$)
and for $t\ne 0$, $J_{i|P_{x_i}\timess\{t\}}
\cong \OO_{P_{x_i}}(1)$. It is easy to see that any surjective homomorphism
$$f_i:\Bl\OO_{P_{x_i}}\Br\uplu{r}\lra \OO_{P_{x_i}}\oplus \CC_{p_i}
\tag 1.39
$$
can be extended to a (surjective) homomorphism
$$F_i:\bigl(\OO_{P_{x_i}\timess U}\bigr)\uplu{r}\lra J_i.
$$
In generally, we can extend any surjective homomorphism
$$f:\pi\sta(L\uplu{r})\lra \oplu{d}{i=1}\pi\sta(L\uplu{r})_{|P_{x_i}}
\lra \oplu{d}{i=1}\Bl \OO_{P_{x_i}}\oplus \CC_{p_i}\Br,\tag 1.40
$$
to a (surjective) homomorphism
$$F: \pi\sta(L\uplu{r})\otimes_{\OO_{\Del}}\OO_{\Del\timess U}\lra
\oplu{d}{i=1} J_i.
$$
Let $V_t=\ker\{F_{|\Del\timess\{t\}}
\}$. Then $V_t$\ is a flat family of torsion free sheaves on $\Del$\
parameterized by $U$. Assuming for general $\sigma_v\mh
\pi\sta(L\uplu{r})\to \oplus^{d}_{i=1}\OO_{P_{x_i}}(1)$, $\Hom
(E_v,E_v(D))^0\ne 0$, then by semicontinuity theorem,
$\Hom(V_t,V_t(D))^0\ne\{0\}$\ for $t\ne0$\ and consequently,
$\Hom(V_0,V_0(D))^0\ne\{0\}$.

Now we seek to find a contradiction by choosing $V_0$\ (i.e. $f$\
in (1.40)) carefully. We first divide the set $\{x_1,\cdots,x_d\}$\
into $2r$\ subsets, say $\Lambda_1,\cdots,\Lambda_{2r}$, such that
each contains either $[d/2r]$\ or $[d/2r]+1$\ points.
We write $f_i=f_i^1\oplus f_i^2$\ according to (1.39).
For $x_i\in \Lambda_{2k-1}$,
we define $f_i^1$\ to be the composition
$$f_i^1:\pi\sta(L^{\oplus r})\mapright{\text{rest.}}
\pi\sta(L^{\oplus r})_{|P_{x_i}}\mapright{\pr_k}\pi\sta(L)_{|P_{x_i}},
\tag 1.41
$$
where $\pr_k$\ is the projection onto the $k^{\text{th}}$\ component
and define $f_i^2$\ to be the composition
$$f_i^2: \pi\sta(L\upr)\mapright{\text{rest}}\pi\sta(L\upr)_{|P_{x_i}}
\mapright{\pr_{k+1}}\pi\sta(L)_{|P_{x_i}}\mapright{\text{ev}} \CC_{p_i},
$$
where $\text{ev}\mh \pi\sta(L)_{|P_{x_i}}\to \CC_{p_i}$\ is the
evaluation map. (Here we agree $\pr_{r+1}=\pr_1$.) For $i\in \Lambda_{2k}$,
we define $f_i^1$\ as in (1.41) while we let $f_i^2$\ to be
$$f_i^2: \pi\sta(L\upr)\mapright{\text{rest}}\pi\sta(L\upr)_{|P_{x_i}}
\mapright{\pr_{k+1}\oplus \pr_{k+2}}\pi\sta(L)_{|P_{x_i}}\oplus
\pi\sta(L)_{|P_{x_i}}\mapright{\text{ev}+\text{ev}} \CC_{x_i}.
$$
($\pr_{r+2}=\pr_2$.) We claim that when $d$\ is sufficiently large, the
sheaf $E\sub\pi\sta(L\upr)$\ that is the kernel of $\oplus_{i=1}^d(f_i^1
\oplus f_i^2)$\ has $\Hom(E,E(D))^0=0$. Indeed, let
$$\tilL_k=L(- \sum_{i\in \Lambda_{2k-1}\cup \Lambda_{2k}} x_i)
$$
be the line bundle on $\Sig$\ of
degree between $l_0-2$\ and $l_0+1$\ and let $L_k=\pi\sta\tilL_k$.
Then, $E$\ is a subsheaf of $\oplus_{i=1}^r
L_h$\ with cokernel $\oplus_{i=1}^d \CC_{p_i}$. Let $s\in\Hom(E,E(D))$.
Then $s$\ induces a homomorphism
$$(s_{ij})_{r\times r}: \oplu{r}{h=1}L_h\lra \oplu{r}{h=1} l_h(D)
$$
with $s_{ij}\in H^0(L_i\upmo\otimes L_j(D))$. Since $L_i\upmo\otimes
L_j$\ is a pull back of line bundle on $\Sig$\ that has
degree $-1$, 0 or 1,
$h^0(L_i\upmo\otimes L_j(D))$\ is bounded by a constant $C_6$\
independent of $d$. On the other hand, by our construction of $E$,
when $i\in\Lambda_{2k-1}$, the composition
$$\oplu{r}{h=1} L_h\mapright{(s_{\ast\ast})}\oplu{r}{h=1} L_h\mapright{
\pr_{k+1}} L_{k+1}\mapright{\text{ev}} \CC_{p_i}
$$
is trivial. Hence for $j\ne k+1$, $s_{jk+1}$\ vanishes on $p_i$\ for all $i\in
\Lambda_{2k-1}$. Now we let $N=2r([C_6]+3)$\ and assume
$d\geq N$. Because $p_i$\ are general and
$$\#(\Lambda_{2k-1})\geq [d/2r]>C_6+2> h^0(L_j\upmo\otimes
L_{k+1})+1,\tag 1.42
$$
$s_{jk+1}$\ must be 0 for $j\ne k+1$.

It remains to show that $s=g_0\cdot\text{id}$\ for some $g_0\in H^0(\OO(D))$.
Let $g_j\in H^0(\OO(D))$\ be sections so that $s_{jj}=g_j\cdot\text{id}\mh
L_j\to L_j(D)$. Let $i\in\Lambda_{2k}$\ and let
$$v_i\ne 0\in\ker\bigl\{\{L_{k+1}\oplus L_{k+2}\}_{|p_i}
\mapright{(\text{ev},\text{ev})} \CC_{p_i}\bigr\}.
$$
Then because $(s_{ij})=\text{diag}\{s_{11},
\cdots,s_{rr}\}$\ is induced from $s\in\Hom (E,E(D))$, we must have
$$(\text{ev},\text{ev})\circ(\pr_{k+1}\oplus \pr_{k+2})\circ(s_{\ast\ast})
v_i=0.
$$
It is straight forward to check that this is equivalent to
$(g_{k+1}-g_{k+2})(p_i)=0$.
Hence, because $p_i$\ are general and $\#(\Lambda_{2k})
> h^0(\OO(D))+1$, we must have $g_{k+1}=g_{k+2}$.
Therefore, $\Hom(E,E(D))^0
=0$. This completes the proof of theorem for $m=0$.

Now we use induction on $m$\ to establish the remaining cases.
The strategy is as follows: We first fix a section
$\Sigp\sub\Del$\ of $\pi\mh\Del\to\Sig$\
of positive self-intersection $\delta$.
Let $E\in\Wr(D,1_m)$\ be any sheaf. We choose a
quotient sheaf $E_{|\Sigp}\to L_E$\ with $L_E$\ a locally
free sheaf of $\OO_{\Sigp}$-modules and
define $\te=\ker\{ E\to L_E\}$. $\te$\ is locally free
with Chern classes
$$I\pri=\det(\te)=I(-\sigp),\quad r_0=\text{rank} L;\tag 1.43
$$
$$d\pri=c_2(\te)=d+\deg L_E+{1\over 2}r_0(r_0-1)\delta -r_0(I\cdott\sigp).
\tag 1.44
$$
Moreover, $\te\in\Wrp(D+\sigp)$\ for a constant $e\pri$\ independent of
$L$\ and $d$. Hence by applying the induction hypothesis to $\Wrp(D+\sigp)$, we
get an upper bound of $\mod\{\te\mid E\in \Wr(D)\}$.
Thus if we understand the correspondence $E\to
\te$\ well, we can translate the estimate of $\mod\{\te\mid E\in \Wr(D)\}$\
to estimate of $\mod\Wr(D)$. We now give the details of this argument.

First, we choose $e_0>0$\ so that $h^1(\Sig, F)=0$\ holds
for all semistable vector bundles $F$\
on $\Sig$\ with $\rank(F)\leq r^2$\ and $\mu(F)\geq e_0$.
Put $e_1=r(e_0+\delta)$. There is a decomposition of
$\Wr(D,1_m)$\ according to whether the restriction of an element $E\in
\Wr(D,1_m)$\ to $\sigp$\ is $e_1$-stable or not.
We denote these sets by $W^+$\ and $W^-$\ respectively.
Let $L_0$\ be a line bundle on $\sigp$\ such that
$H^0_{\sigp}(F\dual\otimes L_0)$\ generates $F\dual\otimes L_0$\
for any $e_1$-stable
rank $r$\ vector bundle $F$\ on $\sigp$\ of degree $I\cdott \sigp$.
Then for any $E\in W^+$, we let $L_E=L_0$\ and fix a surjective
homomorphism $\sigma\mh E\to L_E$.
In case $E\in W^-$, we let
$$0=F_0\sub F_1\sub\cdots\sub F_k=E_{|\sigp}
$$
be the Hardar-Narasimhan filtration of $E_{|\sigp}$. That is
$F_{i+1}/F_i$\ are semistable and $\mu(F_i/F_{i-1})>\mu(F_{i+1}/F_i)$.
We let $i_0$\ be the largest integer so that
$$\mu(F_{i_0}/F_{i_0-1})\geq \mu(F_{i_0+1}/F_{i_0})+e_0.
$$
Such $i_0$\ exists because $E_{|\sigp}$\ is not $e_1$-stable and $e_1>re_0$.
Then by our choice of $e_0$,
$$E_{|\sigp}\cong M_E\oplus L_E,\qquad M_E=F_{i_0}\ \text{and}\
L_E=E_{|\sigp}/F_{i_0}.
$$
We choose our quotient sheaf to be $\sigma\mh E\to L_E$.
Note that $L_E$\ is $re_0$-stable and has degree $\displaystyle
\leq {r_0\over r}I\cdot\sigp$.

Now let $\te$\ be the kernel of $E\to L_E$. Then $\te$\ is locally
free whose first and second Chern classes are given in (1.43) and (1.44).
It can easily be checked that $\te$\ is $e\pri$-stable, $e\pri=e_1
+r(H\cdot \sigp)$, and $\Hom_{\Del}(\te,\te(D+\sigp))^0\ne\{0\}$.
Therefore, we have obtained a map
$$\Psi:\Rr(D,1_m)\lra \bigcup_{d\pri,I\pri}\Wrp(D+\sigp),\tag 1.45
$$
where $d\pri$\ can be any integer and $I\pri$\
can possibly be $I(-\sigp),\cdots,I(-(r-1)\sigp)$.
We wish to find an upper bound on
$$\mod \Psi^{-1}\Bl\Wrp(D+\sigp)\Br
$$
that is independent of
$(d\pri,I\pri)$. We begin with an estimate of $\mod\Psi^{-1}
(\Psi(E))$\ for any $E\in \Rr(D,1_m)$. Because $E$\ belongs to
the exact sequence
$$0\lra F\lra E\lra L\lra 0\tag 1.46
$$
($L=L_E$\ as before) for $M=E_{|\Sigp}/L$, $F_{|\sigp}$\ fits into the exact
sequence
$$0\lra L(-\sigp)\lra F_{|\sigp}\lra M\lra 0.\tag 1.47
$$
On the other hand, elements of $\Psi^{-1}(F)$\ are parameterized by a
subset of $\PP\Ext^1_{\Del}(L,F)$. Since $F$\ is locally free,
$$\align
\dim \Ext^1_{\Del}(L,F)&=\dim \Ext^1_{\Del}(F,L\otimes K_{\Del})
=\dim H^1_{\sigp}({\Cal H}om(F,L)\otimes K_{\Del})\\
&=\dim H^0_{\sigp}(L\dual\otimes F_{|\sigp}(\sigp))
\leq h^0(L\dual\otimes L)+ h^0(L\dual\otimes M(\sigp)).\endalign
$$
Here the last inequality follows from (1.47).
Since $L$\ is $re_0$-stable, $h^0(L\dual\otimes L)$\ is bounded from above by
a constant. In case $E\in W^+$, because $E_{|\Sigp}$\ is $e_1$-stable,
$h^0(L\dual\otimes M(\Sigp))$\ is also bounded from above.
Hence for some constant $C_3$\ depending
only on $e_1$, $r$\ and $I$, we have
$$\mod\Psi\upmo(\Psi(E))\leq C_3,\quad \forall E\in W^+.\tag 1.48
$$
When $E\in W^-$, $h^1(
L\dual\otimes M(\sigp))=0$\ because of our choice of $e_0$. Thus
$$\align
h^0(L\dual\otimes M(\sigp))&=\xx(L\dual\otimes M(\sigp))\\
&=-\deg L\cdot \rk (M)+\deg M\cdot \rk (L)+\rk (M)\cdot \rk (L)
(\delta+1-g).
\endalign
$$
Combined with $\deg L+\deg M=I\cdott \sigp$, we get
$$\mod \Psi^{-1}\Bl\Psi(E)\Br\leq -r\deg L+C_4\tag 1.49
$$
for some constant $C_4\geq 0$\ independent of $E\in W^-$.

Now we use the induction hypothesis. Because $\rk L=r_0<r$, for any
$E\in\Rr(D,1_m)$, either the generic fiber type of $\te$\ is $1_{m-r_0}$\
or it is not of the form $(\alp_1^{\oplus r_1},\alp_2^{\oplus r_2})$.
Hence, with $I\pri=I(-r_0\sigp)$, $1\leq r_0\leq r$\
and
$$C_5=C+C_4+2r^2\delta+2r|I\cdot\sigp|,
$$
we can use theorem 1.2 and the
induction hypothesis to conclude that there is an $N_1$\ and a constant $C_6$\
so that when $d\pri\geq N_1$, we have
$$\mod \Wrp(D+\sigp)\leq \eta(r,d\pri,I\pri)-C_5\tag 1.50
$$
and when $d\pri\leq N_1$, we have
$$\mod \Wrp(D+\sigp)\leq \eta(r,d\pri,I\pri)+C_6.\tag 1.51
$$
We claim that when
$$d\geq N=N_1+r\delta +(2+r)|I\cdot\sigp|+C+C_4+C_6,\tag 1.52
$$
then
$$\mod W^-\leq \eta(r,d,I)-C.
$$
We break the estimate into two cases. In case
$d\pri=c_2(\te)\geq N_1$, then by (1.49) and (1.50),
$$\align
\mod\Psi^{-1}&\Bl\Wrp(D+\sigp)\Br\leq\Bl\eta(r,d\pri,I\pri)-C_5\Br
+\Bl-r\deg L+C_4\Br\\
&=\eta(r,d,I)-2r_0 I\cdot\sigp+\Bl-r_0^2(r-1)+\half r_0(r_0-1)\Br
+r\deg L-C_5+C_4\\
&\leq \eta(r,d,I)-C.\\
\endalign
$$
The last inequality holds because $\deg L\leq {r_0\over r}I\cdot\sigp$.
Now assume $d\pri=c_2(\te)<N_1$. Then
$$\align
\mod\Psi^{-1}&\Bl\Wrp(D+\sigp)\Br\leq\Bl\eta(r,d\pri,I\pri)+C_6\Br+
\Bl-r\deg L+C_4\Br\\
&\leq \eta(r,d,I)+r\deg L+2r|I\cdot\sigp|+r^2\delta+C_6+C_4 \leq \eta(r,d,I)-C.
\endalign
$$
Here we have used the fact that $\deg L\leq -2|I\cdot \sigp|-r\delta
-C_6-C_4-C$\ which follows from (1.44) (1.52) and
$d\pri<N_1$. Now we consider $E\in W^+$. Since
$\mod\Psi\upmo(\Psi(E))\leq C_3$\ from (1.48)
and $c_2(\te)=d+\eta$\ with $\eta$\
a fixed integer independent of $d$, an argument similar to that of $W^-$\
shows that there is an $N\pri$\ such that for $d\geq N\pri$, we have
$\mod W^-\leq \eta(r,d,I)-C$. This establishes the theorem 1.1.
\endpf

\head 2. Degeneration of moduli space
\endhead
\def\rk{\bold{rk}}
\def\rank{\operatorname{rk}}

We now recall briefly the construction of degeneration of moduli and
refer the details of this construction to [GL]. We first fix a very ample
line bundle $H$\ and a line bundle $I$\ on $X$.
Let $C$\ be a Zariski neighborhood of $0\in\spec \CC[t]$. By
choosing a smooth divisor $\Sigma\in |H|$\ we can form a threefold
$Z$\ over $C$\ by blowing up $X\times C$\ along $\Sig\times\{0\}$.
Clearly, $Z_t\cong X$, $t\ne 0$\ and
$Z_0$\ consists of two smooth components $X$\ and a ruled surface $\Del$\
that intersect normally along
$\Sig\sub X$\ and $\sigm\sub \Del$.
For any line bundle $I$\ on $X$\ and integers $r$\ and $d$, let
$\frak M^{r,d}_X$\ be the moduli space of rank $r$\ $H$-semistable
sheaves over $X$\ of $\det E=I$\ and $c_2(E)=d$. Let $\frak M^{r,d}_X
\times C\sta\to C\sta$, $C\sta=C\setminus\{0\}$, be the constant family
over $C\sta$. The degeneration we construct will be a flat family $\Md$\
(over $C$) extending the family $\frak M^{r,d}_X\times C\sta$\
such that the closed points of the special fiber $\Md_0=\Md\times_{C}
\spec \CC[0]$\ are in one-one correspondence with the semistable sheaves
on $Z_0$\ that will be defined shortly

We first introduce the notion of torsion free sheaves on surface $Z_0$:

\definition{Definition 2.1}
A sheaf $E$\ on $Z_0$\ is said to be torsion free at $z\in Z_0$\ if whenever
$f\in\OO_{Z_0,z}$\ is a zero divisor of the $\OO_{Z_0,z}$-modules $E_z$,
then $f$\ is a zero divisor of the $\OO_{Z_0,z}$-modules $\OO_{Z_0,z}$.
The sheaf $E$\ is said to be torsion free if $E$\ is torsion free everywhere.
\enddefinition

Let $E$\ be any coherent sheaf on $Z_0$. We denote by $\eo$\ (resp. $\et$)
the torsion free part of $E_{|X}$ \ (resp. $E_{|\Del}$). We define the
rank of $E$\ to be a pair of integers, $\rk(E)=(\rank(\eo),\rank(\et))$.
When $\rk (E)=(r,r)$, we simply call $E$\ a rank $r$\ sheaf.

Let $\eps\in (0,\half)$\ be a rational number. We define a $\QQ$-ample
divisor $H(\eps)$\ on $Z$\ as follows:
Put $p_X\mh Z\to X$\ be the projection and put
$$H(\eps)=p_X\sta H(-(1-\eps)\Del).
$$
Clearly, for integer $n_0$\ so that $n_0\cdot\eps\in\ZZ$,
$$H(\eps)^{\otimes n_0}=p_X\sta H^{\otimes n_0}(-(n_0-n_0\eps)\Del)
$$
is an ample divisor. In the sequel, we will constantly use the tensor power
$H(\eps)^{\otimes n}$. We agree without further mentioning that in such cases,
$n$\ is always divisible by $n_0$.

Let $\alp=(\alpha_1,\alpha_2)$\ be a pair of rational numbers:
$$\alpha_1=\Bl H(\eps)_{|X}\cdot H(\eps)_{|X}\Br/(H\cdot H),\quad
\alpha_2=\Bl H(\eps)_{|\Del}\cdot H(\eps)_{|\Del}\Br/(H\cdot H).
$$
Note that $\alpha_1+\alpha_2=1$. For any sheaf $E$\ on $Z_0$\ with
$\rk(E)\ne(0,0)$, we define $p_E$\ to be the polynomial
$$p_E={1\over \rk(E)\cdot \alp}\,\xx_E.\tag 2.1
$$
We remark that since $\xx_E(n)=\xx(E\otimes H^{\otimes n})$\
is well-defined for those $n$\ divisible by $n_0$\
and is a restriction of a
polynomial in $n$, we can define $\xx_E$\ to be that polynomial. Once
we have the polynomial $p_E$, we can define the $H(\eps)$-stability
(or $H(\eps)$-semistability) of $E$\ by mimicking the definition 0.3 word by
word.

\definition{Definition 2.2}
A torsion free sheaf $E$\ on $Z_0$\ is said to be $H(\eps)$-stable
(resp. $H(\eps)$-semistable) if whenever $F\sub E$\ is a proper
subsheaf, then $p_F\prec p_E$\ (resp. $\preceq$).
\enddefinition

We fix a line bundle $I$\ on $X$\ and an integer $r\geq 2$.
We let $\xx(n)$\ be the polynomial that depends on $(r,d,I,H,X)$:
$$\xx(n)={r\over 2}n^2(H\cdot H)+n\Bl(H\cdot I)-{r\over 2}(H\cdot K_X)\Br
+(r-1)\xx(\OO_X)+\xx(I)-d.\tag 2.2
$$
$\xx(\cdot)$\ is the Hilbert polynomial of a rank $r$\ sheaf of
$c_1=I$\ and $c_2=d$. We also fix a rational $\eps\in (0,\half)$\
momentarily and the $\QQ$-ample line bundle $H(\eps)$.
For convenience, we will denote by $\E(n)$\ the sheaf
$\E\otimes p_Z\sta H(\eps)^{\otimes n}$\ for any sheaf $\E$\
over $Z_S$. (Here $S$\ is any scheme over $C$\ and $Z_S=Z\times_C
S$.)

We now construct the degeneration $\Md$\ promised at the beginning
of this section. Recall that the moduli space $\frak M^{d,r}_X$\ was
constructed as a GIT quotient of the Grothendieck's Quot-scheme.
Here, we shall adopt the same approach to construct $\Md$.
We first fix a sufficiently large $n$\ and let $\rho=\xx(n)$.
Following A. Grothendieck [Gr], we
define $\underline{\bold{Quot}}^{\xx,\OO^{\rho}}_{Z/C}$\ to be the functor
sending any scheme $S$\ of finite type over $C$\ to the set of all quotient
sheaves $E(n)$\ of $\OO_{Z_S}^{\oplus \rho}$\ on $Z_S$\ flat over $S$\ so that
$\xx_{E_s}(m)=\xx(n+m)$\ for any closed $s\in S$.
$\underline{\bold{Quot}}^{\xx,\OO^{\rho}}_{Z/C}$\
is represented by a scheme $\bold{Quot}^{\xx,\OO^{\rho}}_{Z/C}$\
that is projective over $C$, called Grothendieck's Quot-scheme.
Similarly, we have Grothendieck's Quot-scheme $\bold{Quot}^{\xx,\OO^{\rho}}_X$\
parameterizing all quotient sheaves $\OO_X\uplu{\rho}\to E(n)$\ on $X$\ with
$\xx_E\equiv\xx(\cdot+n)$. Let $\U^{I,d}_X\sub \bold{Quot}^{\xx,\OO^{\rho}}_X$\
be the subset of all $H$-semistable quotient sheaves $E(n)$\ obeying
one further restriction: $\det E=I$. $\U^{I,d}_X$\ is locally closed.
We define $\bold{Quot}^{\xx,\OO^{\rho},I}_{Z/C,H(\eps)}$\ to be the closure of
$\U^{I,d}_X\times C\sta \sub \bold{Quot}^{\xx,\OO^{\rho}}_{Z/C}$\
endowed with the reduced scheme structure and denote by
$\widetilde{\bold{Quot}}^{\xx,\OO^{\rho},I}_{Z/C,H(\eps)}$\
the normalization of $\bold{Quot}^{\xx,\OO^{\rho},I}_{Z/C,H(\eps)}$.
$\widetilde{\bold{Quot}}^{\xx,\OO^{\rho},I}_{Z/C,H(\eps)}$\
has the property that it is normal, projective and flat over $C$.
Finally, we define
$\widetilde{\bold{Quot}}^{\xx,\OO^{\rho},I,ss}_{Z/C,H(\eps)}\sub
\widetilde{\bold{Quot}}^{\xx,\OO^{\rho},I}_{Z/C,H(\eps)}$\ to be the subset of
all closed points whose associated quotient sheaves are $H(\eps)$-semistable.

Clearly, $\widetilde{\bold{Quot}}^{\xx,\OO^{\rho},I,ss}_{Z/C,H(\eps)}$\ depends
on
the choice of $(r,d,n,I,H,\eps)$. In the sequel, $r$, $I$\ and $H$\
will be fixed once and for all. Of course, $d$\ should be viewed as a
variable. For technical reasons, the choice of $\eps$\ will depend on $d$.
After this, we will choose $n$\ sufficiently large (the exact value of $n$\
is irrelevant to our discussion as long as it meets the
requirement   of [Gi, cor.1.3][GL, cor.1.11]).
If all of these are understood, then we will abbreviate
$\widetilde{\bold{Quot}}^{\xx,\OO^{\rho},I,ss}_{Z/C,H(\eps)}$\
to $\U^{d,\eps}$. By abuse of notation, we will call $\E$\ the
universal family of $\U^{d,\eps}$, where $\E(n)$\ is the pullback of the
universal
quotient family on $Z\timec \bold{Quot}^{\xx,\OO^{\rho}}_{Z/C}$.

Let $\SL_C=SL(\rho,\CC)\otimes_{\CC}C$\ be the special linear group
scheme over $C$. Clearly, $\bold{Quot}^{\xx,\OO^{\rho}}_{Z/C}$\ is an
$\SL_C$-scheme.
By our construction, this action lifts to $\U^{d,\eps}$.
Further, we have

\def\Md{\frak M^{d,\eps}}

\pro{Theorem 2.3}
\text{([GL, thm 2.10, 2.11])}\
The good quotient $\Md=\U^{d,\eps}/\!/\SL_C$\ exists. $\Md$\ is
normal, projective
and flat over $C$. Further, for any closed $t\ne 0$,
$\Md_t$\ is isomorphic to the normalization of the moduli scheme
$\MM_X^{r,d}$.
\endpro

To make use of this degeneration, we need to analyze the closed
points of $\Md_0$. Since $\Md_0$\ is a GIT quotient of ${\Cal U}^{d,
\eps}_0$, each point of $\Md_0$\ associates to an equivalent class
of sheaves  $E(n)$\ in ${\Cal U}_0^{d,\eps}$. In the following,
we will find bounds
on $c_1(\eo)$, $c_1(\et)$\ and
$c_2(\et)$\ that are independent of the choice of $\eps$\ and
$E(n)\in \U^{d,\eps}_0$. First we study
$c_1(\eo)$\ and $c_1(\et)$. Following [GL, lem 1.6],
there is a sheaf of $\OO_{\Sig}$-modules $\ez$\ (of rand $r_0$)
such that $E$\ belongs to the exact sequence
$$0\lra E\lra \eo\oplus\et \lra \iota_{\ast}E^{(0)}\lra 0,\tag 2.3
$$
where $\iota\mh \Sig\hookrightarrow X$. Because ${\Cal U}^{d,\eps}$\
is flat over $C$, there are integers $a_1$, $a_2$\ with
$$a_1+a_2=r_0-r\tag 2.4
$$
such that $\det\eo=I(a_1H)$, $\det\et=I_{0|\Del}(a_2\sigm)$, where $I_0=
\ps_X I_{|Z_0}$\ [GL, \S 4].
Then since both $\eo$\ and $\et$\ are quotient sheaves of $E$,
by the $H(\eps)$-stability of $E$, we have
$$a_1\geq (1-\eps){H\cdot(rK_X-2I)\over 2(H\cdot H)};\quad
a_2\geq-(1-\eps){H\cdot(rK_X-2I)\over 2(H\cdot H)}-r.\eqno
$$
Since $H$\ is a very ample divisor on $X$, we may and will assume that
$r$\ divides $H\cdot I$\ and
$$(H\cdot H)\geq 18|K_X\cdot H|+18|I\cdot H|.\tag 2.5
$$
Therefore, $r\geq a_1\geq 0$\ and $0\geq a_2\geq -r$.

The bound of $c_2(\et)$\ is achieved by applying Bogomolov's argument which
shows that when $V$\ is an $H$-stable vector bundle on $X$, then the
restriction of $V$\ to a high degree hyperplane curve is semistable.
We follow the argument in [GL, \S4] and indicate the necessary changes needed
in higher rank case.

\pro{Lemma 2.4}
\text{(cf. [GL, lem 4.3])}
\ There is a constant $e_1$\ independent of $d$\ and $\eps$\ such that the
sheaf
$\eo$\ (on $X$) is $e_1$-stable.
\endpro

\pro{Lemma 2.5}
For any constant $e_1$\ and integer $r$, there is a constant $C_1$\ such
that whenever $V$\ is an $e_1$-stable torsion free sheaf of
rank $\leq r$, then
$$c_2(V)-{\rank(V)-1\over 2\,\rank(V)}c_1(V)^2\geq C_1.\tag 2.6
$$
\endpro

\Proof
The lemma 2.4 is true because any quotient sheaf $Q$\ of $\eo$\ is also a
quotient
sheaf of $E$. Hence the degree of $Q$\ has to satisfy an
inequality, which combined with (2.4) gives us the desired inequality.
The details of the argument can be found in [GL, lemma4.3].
Now we prove lemma 2.5 following the suggestion of the referee.
By Riemann-Roch,
$$\xx(V,V)=\xx(\Ext^{\cdot}(V,V))=
-2r\bigl(c_2(V)-{r-1\over 2r}c_1(V)^2\bigr)+r^2\xx(\OO_X).
$$
Since $V$\ is $e_1$-stable, by lemma 1.8, there is a
constant $C_1$\ such that $\dim \Hom(V,V)$\ and
$\dim\Hom(E,E\otimes K_X)$\ are bounded from above by $C_1$.
Hence
$$c_2(V)-{r-1\over 2r}c_1(V)^2={1\over 2r}\xx(\OO_X)-{1\over 2r}\xx(V,V)
\geq {1\over 2r}\xx(\OO_X)-2C_1.\eqno\square
$$

\pro{Lemma 2.6}
\text{(cf. [GL, lem 4.4])}\ For any constant $e_1$,
there is a constant $C_2$\ such that whenever $V$\ is an $e_1$-stable, rank
$r$\ vector bundle on $X$\ with $\det V=I(aH)$, $|a|\leq r$\ and
that $Q$\ is an $\OO_{\Sigma}$-modules that is a quotient sheaf of
$V_{|\Sigma}$, then we have $\xx(Q)\geq -c_2(V)+C_2$.
\endpro

\Proof
Let $W$\ be the kernel of $V\to Q$. By Riemann-Roch,
$c_1(W)=I+(a-c)[\Sigma]$\ and
$$c_2(W)=c_2(V)+\xx(Q)+\half c( K_X+c H)\cdott H-c(I+aH)\cdott H.
$$
Thus
$$\xx(Q)\geq -c_2(V)+c_2(W)-2r^2 H^2.
$$
On the other hand, since $V$\ is $e_1$-stable, $W$\ is $(e_1+1)$-stable.
So by lemma 2.5, there is a
constant $C_1$\ so that $c_2(W)\geq C_1$.
This completes the proof of lemma 2.6.
\endpf

\pro{Proposition 2.7}
There is a constant $C_3$\ independent of $\eps$\ and $d$\
such that for any $E(n)\in\U^{d,\eps}_0$,
$$ c_2(\et)\leq d+C_3.
$$
\endpro

\Proof
By (2.3), we have
$$\xx_{\et}(\cdot)=\xx_E(\cdot)+\bigl(\xx_{\ez}(\cdot)-
\xx_{\eo}(\cdot)\bigr).
$$
Hence the proposition follows if we can show that the
constant term of $\xx_{\ez}(\cdot)-\xx_{\eo}(\cdot)$\
is bounded from below but this follows from lemma 2.6.
The details of the proof is given in [GL, prop 4.6].
\endpf

Our next goal is to construct Donaldson's line bundle $\LL$\
on $\Md$\ and to establish the following key property of $\LL$:
Whenever $W_0\sub\Md_0$\ is a dimension $c$\ subvariety such that $[\LL]^c
(W_0)>0$, then
$$\mod\bigl\{E^{(2)}\mid E(n)\in W_0\bigr\}=c. \tag 2.7
$$

We first sketch the construction of $\LL$.
The full account of this construction
appeared in [GL,\S5].
For any integer $h\geq 1$, let $D^h\sub Z$\ be a smooth divisor such that
$\pi\mh D^h\to C$\ is smooth, that $D^h_t=\pi^{-1}(t)\in |hH|$\
for $t\ne 0$\ and that $D_0^h\sub \Del\setminus\Sigma$. We call such
$D^h$\ good divisors in $|hH_C(-h\Del)|$, where $H_C=\ps_X H$.
Since $H$\ is very ample, the set of good divisors in $\hdc$\ is base point
free. Associated to each $D^h$\ we can find an \'etale covering $\ctil\to C$\
such that on $\dhct=\dh\timec \ctil$\ there is a line bundle $\thh$\
satisfying $(\thh_v)^{\otimes 2r}=K_{\dh_v}^{\otimes r}\otimes \ps_X
I_{|\dh_v}^{\otimes(-2)}$\ for all closed $v\in \ctil$, where
$K_{\dh_v}$\ is the canonical divisor of $\dh_v$. We remark that
such $\thh$\ exists because $[D^h_v]\cdot I= H\cdot I$\ is divisible by $r$.

We first construct a line bundle on
$\U^{d,\eps}_{\ctil}=\U^{d,\eps}\timec\ctil$\
as follows: Let $\E(n)$\ be the universal quotient family on $Z\timec
\U^{d,\eps}$.
Since $\E$\ is a family of torsion free sheaves flat over $\U^{d,\eps}$,
$\E$\ admits length two locally free resolution near $\dh$. Thus the
restriction of $\E$\ to $\dh\times_C\U^{d,\eps}$\
(denoted by $\E_{|\dh}$) has a length two locally
free resolution also (see [L1]). Let $p_{12}$\ (resp. $p_{13}$; resp. $p_{23}$)
be the projection from $\dh\timec \U_{\ctil}^{d,\eps}$\ to $\dh_{\ctil}$\
(resp. to $\dh\timec \U^{d,\eps}$; resp. to $\U_{\ctil}^{d,\eps}$).
Note that $p_{23}$\ is smooth. Hence
$$R^{\cdot}p_{23\ast}\Bl \ps_{13}(\E_{|\dh})\otimes \ps_{12}\thh\Br\tag 2.8
$$
is a perfect complex on $\U_{\ctil}^{d,\eps}$\ [KM]. Following [KM],
we can define a determinant line bundle
$$\det\Bl R^{\cdot}p_{23\ast}\Bl \ps_{13}(\E_{|\dh})\otimes \ps_{12}\thh\Br
\Br\tag 2.9
$$
on $\U_{\ctil}^{d,\eps}$\ whose inverse we call $\LL_{\U}(\dh)$. If we choose
another good divisor $D^{h\prime}\in\hdc$\
and form the corresponding line bundle
$\LL_{\U}(D^{h\prime})$\ on $\U_{\ctil\pri}^{d,\eps}$,
then since the set of good divisors
in $\hdc$\ is an irreducible set, for any $v\in\ctil$\ and
$v\pri\in \ctil\pri$\ which
lie over the same closed point $t\in C$, the line bundles
$\LL_{\U}(\dh)|\U_{v}^{d,\eps}$\ and
$\LL_{\U}(D^{h\prime})|\U_{v\pri}^{d,\eps}$\
are algebraic equivalent.

\pro{Remark} \rm
Indeed, more is true. There is a single line bundle $\LL_{\U}(h)$\ on
$\U^{d,\eps}$\ such that the line bundles $\LL_{\U}(\dh)$\ on
$\U_{\ctil}^{d,\eps}$\
are pullback of $\LL_{\U}(h)$\ via $\U_{\ctil}^{d,\eps}\to \U^{d,\eps}$.
\endpro

Our next task is to show that under favorable conditions, these line bundles
descend to line bundles on $\Md$. We need the following result of Kempf:

\pro{Lemma 2.8}
\text{(Descent lemma) [DN, thm 2.3]}
Let $\LL$\ be an $\SL_C$\ line bundle on $\U^{d,\eps}$. $\LL$\ descends
to $\Md$\ if and only if for every closed point $w\in\U^{d,\eps}$\
with closed orbit $\SL_C\cdot\{w\}$, the stabilizer $\text{stab}(w)\sub
\SL_C$\ of $w$\ acts trivially on $\LL_w=\LL\otimes k(w)$.
\endpro

We have

\pro{Proposition 2.9}
There is a function $\kappa\mh \ZZ^+\to (0,\half)$\ of which the
following hold:
For any $d$, there is a large $h$\ such that when $\eps\in(0,\kappa(d))
\cap \QQ$\ and $\dh \in\hdc$\ is a good divisor,
then the line bundle $\LL_{\U}(\dh)$\ (on $\U_{\ctil}^{d,\eps}$)
descends to a line bundle on $\Md_{\ctil}=\Md\timec\ctil$. We denote
the descent by $\LL_{\bold{M}}(\dh)$.
\endpro

\proof
It is straightforward to check that $w=E(n)\in \U_{\ctil}^{d,\eps}$\
(over $t\in\ctil$) has closed orbit if and only if $E$\ splits
into direct sum of stable sheaves $F_1\cdots,F_k$.
Then following [L1, p426], the stabilizer $\text{stab}(w)$\ acts trivially
on $\LL_{{\Cal U}}(\dh)_w$\ if and only if
$${1\over \rank(F_1)}c_1(F_1)\cdot D^h_t=\cdots={1\over \rank(F_k)}
c_1(F_k)\cdot D^h_t.
$$
These identities follow if we can prove

\pro{Proposition 2.10}
There is a function $\kappa\mh\ZZ^+\to(0,\half)$\ and a constant $N$\
of which the following hold: Given $d_0$, there is an $h\geq 1$\ such that for
any
$\eps\in (0,\kappa(d_0))$, whenever $d\leq d_0$\ and that
$\E(n)\in\U^{d,\eps}$\
is an $H(\eps)$-semistable sheaves over $t\in C$, then for generic good
divisor $\dh\in\hdc$, $E_{|D^h_t}$\ is semistable.
\endpro

{\it Completion of the proof of proposition 2.9}:
Assume $E=F_1\oplus\cdots\oplus F_k$. By proposition 2.10, there is a good
divisor $D^{\prime h}\in\hdc$\ such that $E_{|D_t^{\prime h}}$\ is semistable.
Then the value ${1\over \rank(F_i)}c_1(F_i)\cdot D_t^h=
{1\over \rank(F_i)}c_1(F_i)\cdot D_t^{\prime h}$\ are identical
for all $i$.
\endpf

Proposition 2.10 will be proved shortly.

\pro{Remark} \rm
Let $t\ne 0\in C$\ be any closed point. Then the line bundle
$\LL_{\U}(\dh)_t$\ on $\U^{d,\eps}_t$\ descends regardless of the choice of
$d$\ and $\eps$\ [L1, p426]. In particular, $\LL_{\text{M}}(\dh)_t$\
always exists on $\Md_t$.
\endpro

Now we explain how to construct global sections of $\LL_{\text{M}}(\dh)_v^{
\otimes m}$\ on $\Md_v$, $v\in\ctil$.
All we need to know about the line bundle $\LL_{\text{M}}(\dh)$\ is
how to calculate its intersection numbers on various subvariety of $\Md$.
So in the following,
we will not distinguish between the line bundles $\LL_{\text{M}}(\dh)_v$\
and $\LL_{\text{M}}(\dh)_{v\pri}$\ (resp. $\LL_{\U}(\dh)_{v}$\ and
$\LL_{\U}(\dh)_{v\pri}$) when $v$\ and $v\pri\in\ctil$\
lie over the same closed point
$t\in C$. By abuse of notation, we will denote either of them by
$\LL_{\text{M}}(\dh)_t$\ (resp. $\LL_{\U}(\dh)_t$).

For any good $\dh\in\hdc$\ and any closed $t\in C$, let $\U^{d,\eps}_t[\dh_t]
\sub\U^{d,\eps}_t$\ be the open set of all $s\in
\U^{d,\eps}_t$\ such that $\E_{s|D^h_t}$\ is semistable. In the
following, we abbreviate $D=\dh_t$.
By restricting $E(n)\in \U^{d,\eps}_t[D]$\ to $D$, we obtain a morphism
$$\Phi_D: \U^{d,\eps}_t[D]\lra \MM^{r,I}(D),\tag 2.10
$$
where $\MM^{r,I}(D)$\ is the moduli scheme of rank $r$\ semistable vector
bundles $V$\ on $D$\ with $\det V=\ps_XI_{|D}$. If we view
$\MM^{r,I}(D)$\ as an $SL(\rho,\CC)$\ scheme with trivial group action,
the morphism $\Phi_D$\ is $SL(\rho,\CC)$\ equivalent.

\pro{Proposition 2.11
\rom{ [Donaldson]}}
There is an ample line bundle $\LL_{D}$\ on $\MM^{r,I}(D)$\ so that
its pull back under $\Phi_D$\ is canonically isomorphic to the
restriction to $\U^{d,\eps}_t[D]$\ of
$\LL_{\U}$.
Further, this isomorphism is $SL(\rho,\CC)$\ equivariant.
\endpro

\Proof
For the details of the proof, the readers are advised to look at [L2, p31].
Though the
author only treated the case $r=2$\ in the proof, the proof of higher
rank case is similar.\endpf

Now let $m$\ be a large positive integer. Since the isomorphism
$$\Phi\sta_D(\LL_D)\cong\LL_{{\Cal U}}|_{\U^{d,\eps}_t[D]}\tag 2.11
$$
is ${ SL}(\rho)$-equivalent, for any
$\xi\in H^0(\MM^{r,I}(D), \LL_D^{\otimes m})$, $\Phi\sta_D(\xi)$\ is an
$SL(\rho)$-invariant section of $\LL_{\U}(D)_t^{\otimes m}$\ on
$\U^{d,\eps}_t[D]$.

\pro{Lemma 2.12}
Let $\dh\in\hdc$\ be any good divisor and for any $t\in C$\ with $D=D_t^h$,
let $\xi\in H^0(\MM^{r,I}(D), \LL_D^{\otimes m})$\ be
any section. Then the pullback section $\Phi\sta_D(\xi)$\
\text (on $\U^{d,\eps}_t[D]$\text )
extends canonically over $\U^{d,\eps}_t$\ to an $SL(\rho,\CC)$-invariant
section. We shall denote this extension (and its descent to
$\Md_t$\ if no confusion is possible) by $\Phi\sta_D(\xi)_{ex}$.
Furthermore,
$$\Phi\sta_D(\xi)_{ex}^{-1}(0)=\bigl(\U^{d,\eps}_t\setminus\U^{d,\eps}_t[D]
\bigr)\cup\bigl\{ F(n)\in \U^{d,\eps}_t[D]\mid \xi(F_{|D})=0\bigr\}.
\tag 2.12
$$
\endpro

\proof
In case $\U^{d,\eps}_t$\ is normal, we can apply
[GL, lemma 5.6][GL, prop. 5.7] and [L2, lemma 4.10] to our situation.
In general, we need to use GIT to prove this lemma
[L1, p435].
\endpf

In the following, we seek to estimate the self-intersection numbers of
$\LL_{\bold{M}}(\dh)$\ on subvarieties
$W\sub\Md_t$\ and to relate the
non-vanishing of such numbers to the estimate of the numbers (2.7).
Our immediate goal is to prove the

\pro{Proposition 2.13}
Let $t\ne0\in C$\ be any closed point and let $W_t\sub\Md_t$\ be an irreducible
variety of dimension $c$, then for sufficiently large $h$\
and for any good $\dh\in\hdc$,
$$[\LL_{\bold{M}}(\dh)]^c(W_t)\geq0.\tag 2.13
$$
Further, if we assume that the general points of $W_t$\
are locally free $H$-$\mu$-stable sheaves, then the strict
inequality holds.
\endpro

\proof
To prove (2.13), it suffices to find divisors $D_1,\cdots,D_c\in
|hH|$\ and sections $\varphi_1,\cdots,\varphi_c$\ of $\LL_{M}(D^h)_t$\
such that $\cap_{i=1}^c\varphi_i^{-1}(0)$\ is a finite set.
But this is obvious because for sufficiently large $h$,
the restriction of each $E\in\Md_t$\ to general $D\i |hH|$\ is
semistable [prop 2.10]. Now we prove the second part of the proposition.
Let $h$\ be large so that for any locally free $E\in\Md_t$,
$H^1(\endo(E)(-hH))=0$. Then for any $D\in |hH|$\ and any locally free
stable $E_1$, $E_2\in W_t$, $E_{1|D}=E_{2|D}$\ implies $E_1\cong E_2$.
We can also assume that the restriction of any $E\in\Md_t$\
to a general $D\in |hH|$\ is semistable.

Choose $D\in |hH|$\ so that $\U^{d,\eps}_t[D]\cap W_t$\ is non-empty. Then
because the line bundle $\LL_D$\ is ample on $\MM^{r,I}(D)$\ and because
$$\Psi_D: \Md_t[D]\cap W_t\lra \MM^{r,I}(D)
$$
$(\Md_t[D]$\ is the image of ${\Cal U}^{d,\eps}_t[D]$\ under the
projection)
is generically one to one, there is a section $\xi\in
H^0(\MM^{r,I}(D), \LL_D^{\otimes m})$, $m$\ large, such that the
extension of the pullback section $\Phi\sta_D(\xi)_{ex}$\ (over
$\Md_t$) is non-trivial over $W_t$\ and
$$\dim\Bigl(\Phi\sta_D(\xi)_{ex}^{-1}(0)\cap W_t\Bigr)=\dim W_t-1.
$$
Since being locally free and stable are open conditions, we can assume
that general points of at least one irreducible component
of $\Phi\sta_D(\xi)_{ex}^{-1}(0)\cap W_t$\ are still locally free and
$H$-$\mu$-stable. Therefore, we can use induction on $\dim W_t$\ to
conclude that for any irreducible component $W_t\pri$\ of
$\Phi_D\sta(\xi)^{-1}_{ex}(0)\cap W_t$, $[\LL_{M}(D^h)]^{c-1}(W_t
\pri)\geq 0$\ and for at least one of these component, this number
is positive. Therefore, the strict inequality (2.13) holds.
\endpf

The converse to the proposition is that if a set $W_t\sub\Md_t$\ with
$\dim W_t=c$\
has the property that
$$[\LL_{\bold{M}}(\dh)]^c(W_t)>0,
$$
then $\mod(W)\geq c$. But this is a tautology since $\Md_t$\ is the
normalization of the moduli scheme.
What we need is a similar result in $t=0$. We will prove

\pro{Proposition 2.14}
Let $W_0\sub\Md_0$\ be any (complete) subvariety of dimension $c$.
Assume for some large $h$\ (given by proposition 2.10) and good
$\dh\in\hdc$\ we have
$$[\LL_{\bold{M}}(\dh)]^c(W_0)>0,
$$
then $\mod\{\et\mid E(n)\in W_0\}=c$.
\endpro

\proof
We prove it by contradiction. Assume $\mod\{\et\mid E(n)\in W_0\}<c$.
Then $\{\et\mid E(n)\in W_0\}$\ can be parameterized by finite
irreducible varieties of dimension at most $c-1$. Let them be
$S_1,\cdots,S_k$\ and let $\E_{1},\cdots,\E_{k}$\ be the corresponding
families. Thanks to proposition 2.10, there is a large $h$\ such that
for any $F\in\{\et\mid E(n)\in W_0\}$, $F_{|D^h_0}$\ is semistable for generic
$\dh\in\hdc$. We fix such an $h$. We choose a $\dh\in\hdc$\
so that $\E_{i,s_i|\dh_0}$\ are semistable for some closed $s_i\in S_i$,
$i=1,\cdots,k$. Since $\LL_{D^h_0}$\ is ample, we can further choose
$\xi\in H^0(\MM^{r,I}(D_0^h), \LL_{D^h_0}^{\otimes m})$, $m\geq 1$,
so that $\xi(\E_{i,s_i|\dh_0})\ne0$\ for all $i$.

Let $\Psi_{D^h_0}\sta(\xi)_{ex}$\ be the extension of the pullback of
$\xi$\ in $H^0(\Md_0, \LL_{\bold{M}}(\dh)_0^{\otimes m})$.
Put $W_0\pri=W_0\cap \Psi_{D^h_0}\sta(\xi)_{ex}^{-1}(0)$.
By our construction, $\dim W_0\pri\leq\dim W_0-1$\ and
$$\mod\{\et\mid E(n)\in W_0\pri\}\leq \max_{i=1,\cdots,k}\{
\dim S_i-1\}\leq \mod \{\et\mid E(n)\in W_0\}-1.
$$
Finally, because
$[\LL_{\bold{M}}(\dh)]^{c-1}(W_0\pri)=m[\LL_{\bold{M}}(D^h)_0]^c(W_0)
>0$. By the induction hypothesis, we have
$\mod\{\et\mid E(n)\in W_0\pri\}\geq c-1$. Therefore,
$$\mod\{\et\mid E(n)\in W_0\}\geq \mod\{\et\mid E(n)\in W_0\pri\}+1\geq c.
$$
The proposition follows because $\mod\{\et\mid E(n)\in W_0\}\leq\dim
W_0= c$.
\endpf

In the remainder of this section, we will give the proof of
proposition 2.10 that is parallel to the treatment
for the rank two situation given in [GL, 5.13].
Let $E(n)\in\U^{d,\eps}_t$\
be any $\he$-semistable sheaf over $t\in C$. When $t\ne0$, then $E$\
is an $H$-semistable sheaf over $X$\ and [MR] tells us that for large
$h$\ and generic $D\in|hH| $, $E_{|D}$\ is semistable. In case $t=0$, namely
when $E$\ is
an $\he$-semistable sheaf on $Z_0$, the situation is quite tricky
because $Z_0$\ is reducible and the divisorial ray $\RR\cdott[D_0^h]$\
is different from $\RR\cdott\he_{|\Del}$. However, it is essential
that $\RR\cdott[D_0^h]$\ and $\RR\cdott\he_{|\Del}$\ become very close
when $\eps$\ becomes small. Before going into the
details of the proof, let us state the following stability criterion of
$\et$.

\pro{Lemma 2.15}
There is a constant $e_2$\ such that for any $d$, $\eps$\
and any $E(n)\in \U^{d,\eps}_0$,
$\et$\ is $\eps e_2$-stable with respect to $\he_{|\Del}$.
\endpro

\Proof
See [GL, 5.14].
\endpf

{\it Proof of proposition 2.10}:
Let $V$\ be the double dual of $\et$. By (2.4) and proposition
2.7, $\det V=I_0(a_2\sigm)$, $-r\leq a_2\leq 0$\ and
$c_2(V)\leq c_2(E)+C_3$, where $C_3$\ is a constant independent of $E$\
and $d$.
Since $I_0\cdot\sigm$\ is divisible by $r$, by tensoring $V$\ with some line
bundle, we can assume $c_1(V)\sim a_2[\sigm]$. Note that $c_2(V)$\
is still bounded by $d_0+C_3$\ possibly with a new constant $C_3$.
Clearly, the proposition will be established if we can show that there is an
$\eps_0$\ and an integer $h$\ such that whenever $\eps<\eps_0$\
and $V$\ is $e_2\eps$-stable with respect to
$\hed$\ as before, then for generic $D\in|h\sigp|$, $V_{|D}$\ is semistable.

The argument we adopt is a direct generalization of Bogomolov's theorem
showing that the restriction of any $\mu$-stable rank two vector bundle
$E$\ to any
smooth hyperplane section of degree $\geq 2c_2(E)+1$\ is stable.
We prove it by contradiction. Assume otherwise. Then there is
a rank $s$\ ($1\leq s\leq r-1$)
quotient vector bundle $Q$\ of $V_{|D}$\ such that
$0=\mu(V_{|D})>\mu(Q)$. Let $W$\ be the kernel of $V\to Q$. Then $W$\
is a locally free sheaf on $\Del$\ with
$c_1(W)\sim a_2[\sigm]-sh[\sigp]$\ and
$$c_2(W)=c_2(V)+\half s(s-1)h^2H^2+\deg Q< c_2(V)+\half s(s-1)h^2H^2.
$$
Thus a simple calculation gives us
$$
2rc_2(W)-(r-1)c_1(W)^2
<2r c_2(V)-\bigl(s(r-s)h^2-(r-1)a_2^2\bigr)H^2.\tag 2.14
$$
Because $c_2(V)\leq d_0+ C_3$, when $\displaystyle
h^2\geq r^2+{2r\over r-1}{d_0+C_3\over
H^2}$\ the right hand side of (2.14) is negative. Therefore, the
Bogomolov's inequality shows that $W$\ is unstable. Let
$$0=W_0\subset W_1\subset \cdots\subset W_n=W\tag 2.15
$$
be the Hardar-Narasimhan filtration of $W$\ such that the sheaves
$F_i=W_i/W_{i-1}$\ are $\mu$-semistable and $\mu(F_i)>\mu(F_{i+1})$.
Let $r_i=\rank (F_i)$\ and let $\Gamma_i$\ be the $\QQ$-divisor supported
on fibers of $\Del\to \Sigma$\ such that
$$c_1(F_i)\sim r_i(b_i\sigm+\Gamma_i).
$$
We let $c_i=\Gamma_i\cdot\sigp/H^2$. Then $b_i$\ and $c_i$\ satisfy
the following inequalities:
$$(e_2+{a_2\over r})\eps\geq \eps b_1+(1-\eps)c_1>\cdots>
\eps b_n+(1-\eps)c_n.\tag 2.16
$$
The first inequality holds because $\et$\ is $e_2\eps$-stable and
the remainder inequalities come from $\mu(F_i)>\mu(F_{i+1})$. On the other
hand, we have $\sum_{i=1}^n c_1(F_i)=c_1(W)$. So
$$\sum_{i=1}^n r_i b_i=a_2-sh,\quad  \sum_{i=1}^n r_ic_i=-sh.\tag 2.17
$$
Finally, we calculate
$$\align
c_2(W)&=\sum_{i<j}c_1(F_i)\cdot c_1(F_j)+\sum_{i=1}^n c_2(F_i)\\
&\geq \half\biggl(\Bigl( \sum_{i=1}^n c_1(F_i)\Bigr)^2- \sum_{i=1}^n
c_1(F_i)^2\biggr)+
\sum_{i=1}^n {r_i-1\over 2r_i} c_1(F_i)^2\\
&=\half(s^2h^2-a_2^2)H^2+\sum_{i=1}^n {r_i\over 2}(b_i^2-2b_ic_i)H^2.
\tag 2.18
\endalign
$$
Here we have used Bogomolov's inequalities $2r_ic_2(F_i)-(r_i-1)c_1(F_i)^2
\geq 0$. Combining (2.18) with $c_2(W)\leq c_2(V)$, we have
$$ (s^2h^2-a_2^2)+\sum_{i=1}^n \BBl r_i(b_i-c_i)^2-r_ic_i^2\BBr \leq
{2(d_0+C_3)\over H^2}. \tag 2.19
$$
In the following, we will argue that there are $h$\ and $\eps_0$\ so that
whenever $0<\eps<\eps_0$, then the only tuples $(b_i,c_i)$\ that
satisfy (2.16)-(2.19) must have $c_i=0$\ for $i=1,\cdots,n-1$.
First of all, let $\Lambda$\ be the set of indices $i$\ so that $c_i>0$.
Then for those $i\in\Lambda$, $c_i\geq 1/rH^2$\ and by (2.16), for
small $\eps$, we have
$$b_i-c_i\leq (e_2+{a_2\over r})+{1\over \eps}(-c_i)<{1\over 2\eps}(-c_i).
\tag 2.20
$$
Thus
$$\sum_{i=1}^n \Bl r_i(b_i-c_i)^2-r_ic_i^2\Br\geq
\sum_{i\in\Lambda} r_i\BBl {1\over 4\eps^2}-1\BBr c_i^2-
\sum_{i\not\in\Lambda} r_i c_i^2.
$$
On the other hand, since $\sum_{i\not\in\Lambda} r_ic_r=
-(sh+\sum_{i\in\Lambda}r_ic_i)$\ and $c_i\leq 0$\ for $i\not\in\Lambda$,
$\sum_{i\not\in\Lambda}r_ic_i^2$\ is bounded from above by
$(sh+\sum_{i\in\Lambda}{r_i}c_i)^2$\ which in turns is no more
than $2s^2h^2+2(\sum_{i\in\Lambda}r_ic_i)^2$. Combined with (2.19),
we must have
$$\sum_{i\in\Lambda}r_i({1\over 4\eps^2}-1)c_i^2-\BBl 2s^2h^2+
2\Bl\sum_{i\in\Lambda}r_ic_i\Br^2\BBr + (s^2h^2-a_2^2)\leq {2(d_0+C_3)\over
H^2}.
\tag 2.21
$$
(2.21) is impossible if we assume
$${1\over 4\eps^2}\geq r^2(r^2h^2+r^2)\cdot H^2+ 2r(d_0+C_3)+4.\tag 2.22
$$
Thus under the assumption (2.22), we must have $c_i\leq 0$\ for
all $i$.

It remains to show that we can choose $h$\ large enough so that
$c_1=\cdots=c_{n-1}=0$. Suppose there are $c_{i_0}<c_{i_1}<0$. Then
$\sum_{i\ne i_0} r_ic_i=-sh-r_{i_0}c_{i_0}$. Again since $c_i\leq 0$, we have
$$\sum_{i=1}^n r_ic_i^2=r_{i_0}c_{i_0}^2+
\sum_{i\ne i_0} r_ic_i^2 \leq r_{i_0}c_{i_0}^2+ (sh+r_{i_0}c_{i_0})^2.
$$
Therefore, from (2.19), we have
$$\align
{2(d_0+C_3)\over H^2}&\geq (s^2h^2-a_2^2)+
\sum_{i=1}^n r_i(b_i-c_i)^2-\sum_{i=1}^n r_ic_i^2\\
&\geq (s^2h^2-a_2^2)-\Bl(sh+r_{i_0}c_{i_0})^2+r_{i_0}c_{i_0}^2\Br
\tag 2.23\\
&\geq (s^2h^2-a_2^2)-\BBl(sh-{1\over r H^2})^2+({1\over rH^2})^2\BBr
= {2sh\over rH^2}-(a_2^2+{2\over r^2 H^2}).
\endalign
$$
Clearly (2.23) is impossible if we choose
$$h\geq 2r(d_0+C_3+r^2H^2)+5.\tag 2.24
$$
Now, we can choose $h$\ large according to (2.24) and then choose
$\eps_0$\ small so
that $\eps_0\leq 1/ 2r$\ and (2.22) holds with $\eps$\
replaced by $\eps_0$. Thus by our previous argument, if $V_{|D}$\
is not semistable, then in the filtration (2.15), all but one
$c_1(W_i/W_{i-1})\cdot [\sigp]=0$. We claim that $c_1(W_n/W_{n-1})\cdot
[\sigp]\ne0$. Indeed, assume $c_j\ne 0$, $j<n$. Then $c_j=-{sh/r_j}$\
and then by (2.19),
$$\BBl \sum_{i\ne j} r_ib_i^2\BBr + r_j(b_j-c_j)^2+ (1-{1\over r_j})s^2h^2
-a_2^2\leq {2(d_0+C_3)\over H^2}.
\tag 2.25
$$
Thus $|b_i|\leq 2\sqrt{d_0+C_3}/\sqrt{H^2}$\ for $i\ne j$\ and $|b_j-c_j|<
2\sqrt{d_0+C_3}/\sqrt{H^2}$.
In particular, we will have
$$\mu(F_j)=(b_j\eps-{(1-\eps)sh\over r_j})H^2\leq b_n\eps H^2=\mu(F_n).
$$
This contradicts to $\mu(F_j)>\mu(F_n)$. Thus we have proved the claim.

The next step is to reconstruct $V$\ from the filtration
$\{W_i\}$. We first
construct a filtration of $V$\ out of the filtration $\{W_i\}$\
by letting ${V}_i \supseteq W_i$\ be the subsheaf of $V$\ so that
$V/{V}_i$\ is torsion free and $\rank(W_i)=\rank({V}_i)$.
We claim that $W_i={V}_i$\ for all $i\leq n-1$. Indeed, let
${V}_i$\ be the first among which ${V}_i\ne W_i$.
Since ${V}_i=W_i$\ on $\Del \setminus D$, we must have
$$c_1({V}_i)=c_1(W_i)+\alp[D],\quad  \alp\geq 1.
$$
On the other hand, $c_1(W_i)=(\sum_{j=1}^i r_jb_j)[\sigm]$\
and $|b_j|\leq 2\sqrt{d_0+C_3}/ \sqrt{H^2}$\ because
of (2.25). Thus
$$\align
\mu({V}_i)&={1\over \rank({V}_i)}c_1({V}_i)
\cdot \hed\\
&= {1\over \rank({V}_i)}\BBl\Bl\sum_{j=1}^i r_jb_j\eps\Br+
{\alp h}\BBr
H^2>\mu(V)+{1\over \rank({V}_i)}e_2\eps\sqrt{H^2},
\endalign
$$
which violates the $e_2\eps$-stability of $V$. Therefore,
${V}_i=W_i$\ for all $i\leq n-1$. In particular, the
filtration
$$0=V_0\subset V_1\subset\cdots\subset V_n=V
$$
has the property that for $i\leq n-1$, $V_i/V_{i-1}$\ are $\mu$-semistable and
$c_1(V_i/V_{i-1})\sim r_ib_i[\sigm]$. Let $F_i=V_i/V_{i-1}$.
We intend to use induction on the rank $r$\ to complete the proof of the
proposition. In order to do this, we need to show that $F_n$\ is
$\mu$-semistable and
$$c_2(V_i/V_{i-1})-{r_i-1\over 2r_i}c_i(V_i/V_{i-1})^2\leq d_0+C_3\tag 2.26
$$
for all $i\leq n$.
We show $F_n$\ is $\mu$-semistable by showing that $r_n=1$. Indeed,
a combination of (2.25) (with $j=n$) and (2.24) guarantees $r_n=1$.
Thus $F_n$\ is stable. Next, we have
$$\align
c_2(V)-{r-1\over 2r}c_1(V)^2&=\sum_{i=1}^nc_2(F_i) +\sum_{i<j} c_1(F_i)\cdott
c_1(F_j)+ (-\half+{1\over 2r})\Bigl(\sum_{i=1}^n c_1(F_i)\Bigr)^2\\
&=\sum_{i=1}^n\BBl c_2(F_i)-{r_i-1\over 2r_i}
c_1(F_i)^2-{1\over 2r_i}c_1(F_i)^2\BBr
+{1\over 2r}\Bl\sum_{i=1}^n c_1(F_i)\Br^2\\
&=\sum_{i=1}^n\Bl c_2(F_i)-{r_i-1\over 2r_i} c_1(F_i)^2\Br+
{1\over 2r}c_1(V)^2
+{1\over 2}\sum_{i=1}^n{1\over r_i}r_i^2 b_i^2H^2.
\endalign
$$
Because each $c_2(F_i)-{r_i-1\over 2r_i} c_1(F_i)^2$\ is non-negative,
(2.26) must be true.
Therefore, we can apply the
induction argument to $V_i/V_{i-1}$\ to conclude that we can find
large $h$\ and small $\eps_0$\ so that for any $\eps<\eps_0$, we
must have $(V_i/V_{i-1})_{|D}$\ semistable for generic $D\in|h\sigp|$.
Since $\deg(V_i/V_{i-1})_{|D}=0$, $V_{|D}$\ must be semistable also.
This completes the proof of proposition 2.10.
\endpf

\head 3. Main theorems
\endhead
\def\erdi{\eta_X(r,d,I)}
\def\md{\frak M^{d,\eps}}

In this section, we will prove our main theorems.
We will show that when the second Chern class $d$\ is large enough, then the
moduli
scheme $\mrdx$\ ($=\mrdx(I,H)$) is smooth at a dense open subset. We shall
further show that $\mrdx$\
is normal and for any constant $C$, there is
an $N$\ depending on $(X,I,H,r,C)$\ such that whenever $d\geq N$, then
$$\codim\Bl\text{Sing}\,\mrdx,\mrdx\Br\geq C.
$$
Finally, we will
investigate the dependence of the moduli scheme $\mrdx(I,H)$\ on
the polarization $H$. In case $r=2$, Qin's work
[Qi] shows that for any two polarizations $H_1$\ and $H_2$, the corresponding
moduli spaces $\mrdx(I,H_1)$\ and $\mrdx(I,H_2)$\ are birational when $d$\
is sufficiently large. Here, we shall demonstrate that the
similar phenomena also occurs in high rank cases.
But first, we shall continue our discussion of the degeneration to finish
the proof our main technical theorem 0.1.

For the moment, we shall
keep the notation developed in \S 2. For any divisor $D\sub X$,
we define $\she\rd$\ be the set of all $e$-stable
(with respect to the fixed $H$)
rank $r$\ sheaves $E$\ of $\det E=I$\ and $c_2(E)=d$\ and define
$$\she\rd(D)=\{E\in\she\rd \mid \Hom(E,E(D))^0\ne\{0\}\}.
$$
Similarly, we define $\vecx\rd$\ and $\vecx\rd(D)$\ to be the subsets of
locally free sheaves in $\she\rd$\ and $\she\rd(D)$\ respectively.
For technical reasons, we will first attack the set $\vecx\rdm(D)$\ which is
the set of $\mu$-stable locally free sheaves $E$\ with the mentioned
constraint on $c_1$, $c_2$\ and $h^0$. Namely, $\vecx\rdm(D)=\vecx^{r,d}_{0,I}
(D)$. We shall prove

\pro{Theorem 3.1}
For any choice of $r$, $I$\ and $D$, and any constant $C_1$,
there is a constant $N$\ such that whenever $d\geq N$, we have
$$\mod \vecx\rdm(D)\leq \eta_X(r,d,I)-C_1.
$$
\endpro

\Proof
Clearly, $\vecx\rdm(D)$\ is a subset of $\mrdx$. Since being
locally free and $\mu$-stable are open conditions and having
non-vanishing $\Hom(E,E(D))^0$\ is a closed condition, $\vecx\rdm(D)$\
is a locally closed subset of $\mrdx$. Let $\bold{A}\sub\mrdx$\ be
the closure of any irreducible component of $\vecx\rdm(D)$.

In the following, we seek to utilize the degeneration
$\md\to C$\ (of the normalization
of $\mrdx$) constructed in theorem 2.3.
When $t\ne 0$, $\md_t$\ is just the normalization of $\mrdx$. For such
$t$, we let $W_t\sub\md_t$\ be the preimage of $\bold{A}\sub\mrdx$.
$\cup_{t\ne0}W_t$\ is a constant family over $C\sta$.
We then let $W$\ be the closure of $\cup_{t\ne0}W_t$\ in $\md$\
and let $W_0$\ be the special fiber of $W$\ over $0\in C$.

Here is our strategy: Take a large $h$\ and a good $\dh\in\hdc$.
By proposition 2.13, for any $t\ne0$\ and $c=\dim W_t$, the top
self-intersection
number $[\LL_{\dh}]^c(W_t)>0$.
Then since $W$\ is flat and proper over $C$,
$[\LL_{\dh}]^c(W_0)>0$.
Therefore, according to proposition 2.14,
$$\mod\{\et|E\in W_0\}=c.\tag 3.1
$$
On the other hand, since every $E\in W_t$, $t\ne0$, has
non-vanishing $\Hom(E,E(D))^0$, the upper-semicontinuity
theorem tells us that there is a divisor $D\pri\sub \Del$\ (we can make it
being independent
of the choice of $W_0$) such that for any $E\in W_0$,
$$\Hom_{\Del}(\et,\et(D\pri))^0\ne\{0\}.\tag 3.2
$$
Thus, by applying theorem 1.1, we get an upper
bound of (3.1) and hence an upper
bound of $c$. We now fill in the details of this approach.

To establish (3.2) for some $D\pri\sub \Del$, we argue as follows:
First of all, let $E\in W_0$\ be any point. Since $W$\ is flat over
$C$, there is a smooth affine curve $S$\ over $C$\ and a flat
family of torsion free sheaves $E_S$\ on $Z_S=Z\!\times_S\!S$\ such
that for any closed $s\in S$\ over $t\ne0\in C$, $E_s\in W_t$\
and further, there is a closed $s_0\in S$\ over $0\in C$\ so that
$E_{s_0}=E$. For any integer $k$, we
consider the divisor $D_C-k\Del$\ on $Z$,
where $D_C=p_X^{-1}(D)$, and the pullback divisor
(of $D_C-k\Del$) on $Z_S$\ which we denote by $D_k$. Clearly, the
restriction of $D_k\otimes k(s_0)$\ to $X\sub Z_0$\ is
$D-kH$. Now consider the vector space $\Hom_{Z_S}(E_S,E_S(D_k))^0$.
By assumption, for general $s\in S$, $\Hom_{Z_s}(E_s,E_s(D_k))^0\ne0$.
Thus $\Hom_{Z_S}(E_S,E_S(D_k))\ne\{0\}$. Let
$w\in \Hom_{Z_S}(E_S,E_S(D_k))^0$\
be a non-trivial section and let $\xi$\ be the uniformizing parameter
of $S$\ at $s_0$. Then because $E_S$\ is flat over $S$,
there is an $n\geq 0$\ such that the restriction
of $w/\xi^{n}$\ to $Z_{s_0}$\ gives rise to a non-trivial
homomorphism $\varphi\mh E_{s_0}\to E_{s_0}(D_k)$.

Next, because $E^{(1)}$\ is a quotient sheaf of $E$,
$\varphi$\ induces a homomorphism $E\to E^{(1)}(D_k)$\ and further because
$E^{(1)}$\ is torsion free, it comes from
$\varphi_1\mh E^{(1)}\to E^{(1)}(D_k)$. Similarly, we have $\varphi_2\mh
E^{(2)}\to E^{(2)}(D_k)$.
Because $E$\ is torsion free, at least one $\varphi_i$\ is non-trivial.
Now we claim that we can choose a $k$\ (independent of $d$\ and $\eps$)
so that $\varphi_1$\
is always trivial. Indeed, we first choose $k$\ so that $H\cdot(D-kH)<0$.
Then since $\det\varphi_1\in H^0(\OO_X(rD-rkH))=\{0\}$, $\det \varphi_1$\
is trivial.
If we let $A\subset \eo$\ be the kernel of $\varphi_1$, then $\eo/A$\
is torsion free and further, there is a $\psi$\ making the following
diagram commutative
$$\CD
A @>>> \eo @>>> \eo/A @>>> 0.\\
@. @VV{\varphi_1}V @V{\psi}VV\\
 @. \eo(D_k) @= \eo(D_k)\\
\endCD
\tag 3.3
$$
On the other hand, by lemma 2.4,
there is a constant $e_1$\ independent of $d$\ and $\eps$\ such that
$\eo$\ is $e_1$-stable. Thus if $\varphi_1\ne0$, then
$0<\rank(\eo/A)<\rank(\eo)$\ and $E^{(1)}/A$\ is both a subsheaf of
$E^{(1)}(D_k)$\ and a quotient sheaf of $E^{(1)}$. Therefore,
$$\mu(\eo(D_k))+{1\over \rank(\eo/A)}\sqrt{H^2}\cdot e_1>
\mu(\eo/A)>
\mu(\eo)-{1\over \rank(\eo/A)}\sqrt{H^2}\cdot e_1.
$$
A straightforward calculation shows that this is impossible if we let
$$k>{1\over H^2}(D\cdot H+2\sqrt{H^2}\cdot e_1).\tag 3.4
$$
Hence $\varphi_1$\ must be trivial.

{}From now on, we fix such a $k$. Then our previous argument shows that all
$V\in\Theta=\{\et\mid E\in W_0\}$\
have non-vanishing $\Hom_{\Del}(V,V(D_k))^0$.
As we explained, our intention is to apply theorem 1.1 to the set $\Theta$\
to get the bound:
$$\mod\Theta\leq \eta_X(r,d,I)-C_1,\quad d\gg0.\tag 3.5
$$
First of all, all $V\in\Theta$\ are $e_2$-stable by lemma
2.15 and have $\det V=I_2$, $I_2\in\Lambda=\{I_0(-r\sigm),\cdots,I_0\}$.
Next, for each $I_2\in\Lambda$, there is an $\eps_0(I_2)>0$\
specified by theorem 1.1. We let $\eps_0=\min_{I_2\in\Lambda}\{\eps_0(I_2)\}$.
Then for any $\eps$\ smaller than the $\eps_0$, the ample
divisor $\hed$\ on $\Del$\ satisfies the condition of theorem 1.1.
In order to apply theorem 1.1, we need to know
that the general element of $\Theta$\ is locally free, which certainly
is quite delicate in general. The solution we propose is to use the double
dual operation to relate any sheaf $F\in\Theta$\ to its double dual
$\F(F)=F^{\vee\vee}$. $F^{\vee\vee}$\
is always locally free because $\Del$\ is a smooth surface.
Assume $d_2=c_2(\F(F))$, then $d_2\leq c_2(F)$\ and the equality holds if
and only if $F$\ is locally free. Following the notation introduced at the
beginning of \S 1, we have
$$\F:\Theta\lra \bigcup_{d_2\in\ZZ\,; I_2\in\Lambda}\Wrt(D_{k|\Del}).
$$
(We use $\frak A^{\cdot}_{\cdot}$\ to denote sets related to $\Del$\
and use ${\Cal V}^{\cdot}_{\cdot}$\ to denote sets related to $X$.)
Here we have used the fact that $\Hom\Bl F,F(D_{k|\Del})\Br^0\ne0$\ implies
$\Hom\Bl F\dual,F\dual(D_{k|\Del})\Br^0\ne0$.
Next, we divide $\Theta$\ into subsets $\Theta_{d_1}$\ according to the
value of the second Chern class of $F\in\Theta$. Then,
$\Theta=\cup \Theta_{d_1}$. We have the following estimate which will
be proved shortly.

\pro{Lemma 3.2}
For any $V\in\Wrt$, $\mod\Bl \F^{-1}(V)\cap \Theta_{d_1}\Br\leq
(r+1)(d_1-d_2)$.
\endpro

Now we are ready to complete the proof of the theorem.
First of all, by applying theorem 1.1 to the set $\Wrt(D_{k|\Del})$, we
know that for any constant $C_2$, there is an $N_2$\ such
that whenever $d_2\geq N_2$, we have
$$\mod\Wrt(D_{k|\Del})\leq \eta_{\Del}(r,d_2,I_2)-C_2, \quad I_2\in\Lambda.
\tag 3.6
$$
To control the left hand side of (3.6) for small $d_2$, we invoke
theorem 1.5 to get
$$\mod\Wrt \leq \eta_{\Del}(r,d_2,I_2)+C_3,\quad I_2\in\Lambda
\tag 3.7
$$
where $C_3$\ is a constant. Another estimate
we need was established in proposition 2.7,
$$c_2(F)\leq d+ C_4, \quad \forall F\in\Theta.\tag 3.8
$$
The proof of (3.5) then goes as follows:
For any constant $C_1$, we let $C_2$\ be such that
$$C_2\geq C+\eta_{\Del}(r,C_4,I_2)-\eta_X(r,0,I),\quad
\forall I_2\in\Lambda\tag 3.9
$$
and let $N_2$\ be the constant that makes (3.6)
hold. We then let $N$\ be so that
$$(r+1)(N-N_2)+\eta_{\Del}(r,N_2,I_2)+C_3\leq \eta_X(r,N,I)-C_1.\tag 3.10
$$
We claim that when $d\geq N$\ and $\eps<\eps_0$, then
(3.5) holds. Indeed, let $d_1\leq d+C_4$\ be any integer. Then
for $d_2\leq N_2$, by (3.7) and (3.10),
$$\align
\mod\biggl(\Theta_{d_1}\bigcap & \F^{-1}\Bl\Wrt(D_{k|\Del})\Br\biggr) \\
 &\leq (r+1)(d_1-d_2) +
\mod\biggl(\F(\Theta_{d_1})\bigcap\Wrt(D_{k|\Del})\biggr) \\
&\leq (r+1)(d_1-d_2)+ \eta_{\Del}(r,d_2,I_2)+C_3\leq
\eta_X(r,d,I)-C_1.\endalign
$$
Assume $d_2\geq N_2$. By (3.6) and (3.9), we have
$$\align
\mod \biggl(\Theta_{d_1}\bigcap & \F^{-1}\Bl\Wrt(D_{k|\Del})\Br\biggr)
\leq \mod\Wrt(D_{k|\Del})+(r+1)(d_1-d_2)\\
&\leq \eta_{\Del}(r,d_2,I_2)-C_2+(r+1)(d_1-d_2)\leq \eta_X(r,d,I)-C_1.\endalign
$$
Thus we have established (3.5).

To finish the proof of the theorem, it suffices to show that
$\mod\Theta=\mod W_t$. For this, we will use the Donaldson's
line bundle $\LL_{\dh}$. First of all, for any $d$, we choose
$\eps<\min(\eps_0,\kappa(d))$. ($\kappa(d)$\ was specified in proposition
2.10.) We then apply proposition 2.13 to the set $W_t$, $t\ne0$. Proposition
2.13 asserts that with $c=\dim W_t$, $[\LL_{\dh}]^c(W_t)>0$. Since $W$\
is flat and proper over $C$, we have
$$[\LL_{\dh}]^c(W_0)=[\LL_{\dh}]^c(W_t)>0.
$$
In particular, proposition 2.14 tells us that then
$$\mod\{\et\mid E\in W_0\}=c.
$$
Therefore, combined with inequality (3.5), we have that for $d\geq N$,
$$\dim\{ E\in\mrdx\mid \Hom_X(E,E(D))^0\ne0\}\leq \eta_X(r,d,I)-C.
$$
This completes the proof of theorem 3.1.
\endpf

Before we go any further, let us finish the proof of lemma 3.2.

{\it Proof of lemma 3.2}:
The situation when $r=2$\ was proved in [L1, p461]. In general, let $E$\ be
any rank $r$\ torsion free sheaf
and let $V=\F(E)$. Then $E$\ is uniquely determined by the quotient
sheaf $V\to V/E$, where $V/E$\ is supported on a discrete
set and of length $\len(V/E)=c_2(E)-c_2(V)$. Therefore,
$\F^{-1}(V)\cap\{\hbox{sheaves of } c_2=c_2(V)+c\}$\
is exactly the set of all quotient sheaves $V\to A$\ such that
$A$\ is supported on a discrete set and $\len (A)=c$. Let $\Quv$\ be the
Grothendieck's Quot-scheme of all quotient sheaves $A$\ of $V$\ with
$\len(A)=c$. $\Quv$\ is projective by [Gr, p13]. Observe also that when $A$\
is supported on $c$\ distinct points, then by [Gr, p21],
$$\dim T_A\Quv=(r+1)c.\tag 3.11
$$
Thus the lemma will be established if we can show that for any quotient
sheaf $A_0\in\Quv$, there is a deformation $A_t$\ of $A_0$\ such that
for generic $t$, $A_t$\ supports on $c$\ distinct points [L1, p461].
In the following,
we will demonstrate how to construct such a deformation.

Clearly, this is a local problem. Let $U$\ be a classical neighborhood
of $0\in \CC^2$\ with coordinate $z=(z_1,z_2)$. Assume $A_0$\ is
a quotient sheaf of $\OO_U^{\oplus r}$\ of length $c$\ supported at the
origin 0. Let $E=\ker\{\OO_U^{\oplus r}\to A_0\}$. Along the
lines of the argument given in [L1, p462], we can show that there are
$f_1,\cdots,f_n\in\OO_U^{\oplus r}$\
such that $\{f_i\}_{i=r+1}^n$\ are divisible by $z_1$\ and $\{f_i\}_1^n$\
generate the submodule $E$.

Next, we define
$$f_i(z,t)=\cases
f_i(z),&1\leq r\leq r;\\
(z_1-t){f_i(z)\over z_1},&r+1\leq i\leq n.\endcases
\tag 3.12
$$
We then define a submodule $E_D\subset \OO_{U\times D}^{\oplus r}$, where
$D$\ is a small disk with parameter $t$, by
$$E_D=(f_1(z,t),\cdots,f_n(z,t))\cdot \OO_{U\times D}\subset \OO_{U\times
D}^{\oplus r}.\tag 3.13
$$
Let $A_D=\OO_{U\times D}^{\oplus r}/E_D$. $E_D$\ and $A_D$\ can be
viewed as families of sheaves parameterized by $D$. It
is easy to see that
when $E_D\otimes k(0)$\ is torsion free, then $A_D\otimes k(0)=A_0$\ and
for $t$\ small, $A_D$\ is a (flat) deformation of $A_0$.
Now we check that $E_D\otimes k(0)$\ is torsion free. Suppose there are
$h\in E_D$\ and $f\in \OO_U$\ such that $f\cdot h=th\pri$\ for some $h\pri
\in E_D$. Let
$$h=\sum_{i=1}^n g_1(z,t)\cdot f_i(z,t).
$$
Then the fact that $f(z)\cdot h\equiv0\,\text{mod}(t)$\ in $\OO_{U\times D}
^{\oplus r}$\ and that $f_1(z)/(z_1),\cdots,f_r(z)/(z_1)$\ generate a rank
$r$\ \ $\OO_{U}/(z_1)$-module implies
$$z_1\mid g_i(z,0),\quad i=1,\cdots,r;\tag 3.14
$$
$$\sum_{i=1}^r{g_i(z,0)\over z_1}\cdot f_i(z)+\sum_{i=r+1}^ng_i(z,0)\cdot
{f_i(z)\over z_1} \equiv 0.\tag 3.15
$$
Further, if we write $g_i(z,t)=\alpha_i(z)+t\beta_i(z,t)$, the following
identities hold in $\OO_{U\times D}^{\oplus r}$:
$$\align
h&=\sum_{i=1}^r\Bigl(\alpha_i(z)f_i(z)+t\beta_i(z,t)f_i(z)\Bigr)\\
&\quad\quad +\sum_{i=r+1}^n\Bigl(\alpha_i(z) (z_1-t){f_i(z)\over z_1}+
t\beta_i(z,t) (z_1-t){f_i(z)\over z_1}\Bigr)\\
&=(z_1-t)\Bigl(\sum_{i=1}^r{\alpha_i(z)\over z_1}f_i(z)+\sum_{i=r+1}^n
\alpha_i(z){f_i(z)\over z_1}\Bigr)+\\
&\quad\quad+t\Bigl(\sum_{i=1}^r\bigl({\alpha_i(z)\over z_1}+\beta_i(z,t)\bigr)
f_i(z)+\sum_{i=r=1}^n\beta_i(z,t)(z_1-t){f_i(z)\over z_1}\Bigr)\\
&=th^{\prime\prime},
\endalign
$$
where $h^{\prime\prime}$\ obviously belongs to $E_D$. Since $E_D$\ is a
submodule of $\OO_{U\times D}^{\oplus r}$, $h$\ must be equal to
$th^{\prime\prime}$\ in $E_D$. Therefore, $f\cdot h_{|t=0}=0$\
implies $h_{|t=0}=0$\ in $E_D\otimes k(0)$\ or that $F_D\otimes k(0)$\
is torsion free.

In general, $A_t$\ is not supported on $c$\
distinct points. But at least we expect that $A_t$\ is simpler than $A_0$,
say supp($A_t$) has at least two distinct points. In the following, we will
show that this is indeed the case.
Without loss of generality, we can assume that $f_i(z)$\ all vanish at
the origin. (Since otherwise, $A_0$\ is essentially a quotient sheaf of
$\OO_U^{\oplus (r-1)}$\ and we can use induction on $r$\ to take care
this situation.)
For small $t$, the equation
$$f_1(t,z_2)\wedge\cdots\wedge f_r(t,z_2)=0
$$
has solutions, say $z_2=w_t$, because
$f_1(0)\wedge\cdots\wedge f_r(0)=0$\
and $f_1(0,z_2)\wedge\cdots\wedge f_r(0,z_2)\ne 0$\ for generic $z_2$.
Note that $(t,w_t)\in \text{Supp}(A_t)$.
If Supp($A_t$) is a single point, then ${f_{r+1}(z)/ z_1},\cdots
{f_n(z)/ z_1}$\ must generate $\OO_U^{\oplus r}$\ at the origin. Thus by
discarding some extra terms, we will have $n=2r$\ and further,
by eliminating terms in $f_1,\cdots,f_r$\ that involves $z_1$\ by using
combination of $f_{r+1}\cdots,f_n$, we can assume $z_2|f_1(z),\cdots,
z_2|f_r(z)$. Therefore, we can consider the deformation of $A_0$\
derived from
$$E_D\pri=\Bigl((z_2-t){f_1(z)\over z_2},\cdots,(z_2-t){f_r(z)\over z_2},
f_{r+1}(z),\cdots,f_n(z)\Bigr).
$$
In case $\text{Supp}(A_t\pri)$\ is still a
single point for generic $t$, then $\Bigl({f_1(z)\over z_2},\cdots.
{f_r(z)\over z_2}\Bigr)$\ will generate $\OO_U^{\oplus r}$\ at 0 also. In
particular, $A_0=\oplus^r \CC$\ and then the desired deformation can be written
by hand.
\endpf

In the remainder of this section, we will complete the proof of the theorems
stated at the beginning of this paper.
We first investigate the sets $\she\rd$\ and $\she\rd(D)$\ introduced at
the beginning of this section. We shall prove

\pro{Theorem 3.3}
For any choice of $r$, $I$\ and $D$\ and any choice of constants $e$\ and $C$,
there is an integer $N$\ such that whenever $d\geq N$, then we have
$$\mod\she\rd=\eta(r,d,I),\tag 3.17
$$
$$\mod\she\rd(D)\leq\erdi-C.\tag 3.18
$$
\endpro

{\it Proof of (3.17)}.
Let $\vecx\rd=\she\rd\cap\{\hbox{locally free sheaves}\}$\
and let $\vecx\rd(D)=\vecx\rd\cap\she\rd(D)$.
Clearly, (3.17) is a stronger statement than
$$\mod\vecx\rd=\erdi, \tag 3.19
$$
which in turn is stronger (in case $e>0$) than
$$\mod\vecx\rdz=\erdi.\tag 3.20
$$
Our strategy is first to prove statement (3.20) and then prove (3.19)
and (3.17). We proceed by induction on the rank $r$. (3.17) and
(3.20) are trivial when $r=1$.
For $r\geq 2$\ and $E\in\vecx\rdz$,
the Kodaira-Spencer-Kuranishi deformation theory tells
us that there is a holomorphic map
$$f: U\subset H^1(X,\endo(E))\lra H^2(X,\endo(E)),
$$
where $U$\ is an (analytic) neighborhood of the origin, such that
$f^{-1}(0)$\ is the versal deformation space of $E$. Since $h^0(\endo(E))=0$\
(since $E$\ is $\mu$-stable),
$$\mod\Bl\vecx\rdz,[E]\Br \geq h^1(\endo(E))-h^2(\endo(E)),
$$
and when $h^2(\endo(E))=0$, $\mod\Bl\vecx\rdz,[E]\Br=h^1(\endo(E))$.
Next, by Riemann-Roch, one calculates $\xx(\endo(E))=\erdi$. Thus one gets
$$\mod\Bl\vecx\rdz,[E]\Br \geq \erdi.\tag 3.21
$$
On the other hand, since $h^2(\endo(E))=h^0(\endo(E)\otimes K_X)$,
by theorem 3.1, there is an $N$\ such that whenever $d\geq N$, we have
$$\mod\{E\in\vecx\rdz\mid h^0(\endo(E)\otimes K_X)>0\}\leq \erdi-1.
$$
Therefore, for generic $E\in\vecx\rdz$, $\mod\Bl\vecx\rdz,[E]\Br=\erdi$.
Thus we have proved (3.20) provided $d\geq N$.
To further attack (3.19) and (3.17), we need
the following estimate which is interesting in its own right.
\def\rdi{^{r_i,d_i}_{e_3,I_i}}

\pro{Theorem 3.4}
For any choice of $r$, $I$\ and two constants $e_1>e_2$, there is a constant
$C\pri$\ such that
$$\mod\Bl\she\rdo\setminus\she\rdt\Br \leq (2r-1)d+C\pri. \tag 3.22
$$
\endpro

\proof
Let $E$\ be any torsion free sheaf in $\she\rdo\setminus\she\rdt$.
Since $E$\ is not
$e_2$-stable, there is a torsion free subsheaf $F_1\sub E$\ such that
$E/F_1$\ is torsion free and that
$\mu(F_1)\geq \mu(E)+e_2\sqrt{H^2}/\rank(F_1)$. Because
$E$\ is $e_1$-stable, $\mu(F_1)$\ is bounded from above by $\mu
(E)+e_1\sqrt{H^2}/\rank(F_1)$. Combined, we get
$${1\over r}I\cdot H+\frac{1}{r_1}e_1\sqrt{H^2}
<{1\over r_i}I_i\cdot H< {1\over r}I\cdot H+
{1\over r_i}e\sqrt{H^2},\tag 3.23
$$
where $r_i=\rank F_i$, $d_i=c_2(F_i)$\ and $I_i=\det F_i$\ with
$F_2=E/F_1$. Note
that $E$\ belongs to the exact sequence
$$0\lra F_1\lra E\lra F_2\lra 0.\tag 3.24
$$
We call $(r_i,d_i,I_i)$\ admissible if they do come from
(3.24) with $E\in\she\rdo\setminus\she\rdt$.
We claim that $F_i$\ are $e_3$-stable with $e_3=e_1+|e_2|$.
Indeed, let $L\subset F_1$\ be any subsheaf.
Because $L$\ is also a subsheaf of $E$,
$$\mu(L)<\mu(E)+{1\over \rank(L)}e_1\sqrt{H^2}\leq\mu(F_1)-{1\over r_1}e_2
\sqrt{H^2}+{1\over \rank(L)}e_1\sqrt{H^2}.
$$
Thus, $F_1$\ is $e_3$-stable. $F_2$\ is
$e_3$-stable for the same reason. Therefore, $F_i\in\she\rdi$. Finally,
because of (3.24),
$$\align
\mod\Bl\she\rdo\setminus\she\rdt\Br\leq \sup_{(r_i,d_i,I_i)}&\bigl\{
\mod \Bl\she^{r_1,d_1}_{e_3,I_1}\Br+\mod \Bl\she^{r_2,d_2}_{e_3,I_2}\Br
\\
&+\sup\{\dim\Ext^1(F_2,F_1)\mid F_i\in\she\rdi\}\bigr\},\tag 3.25
\endalign
$$
where the supremum is taken over all admissible  tuples $(r_i,d_i,I_i)$.
Note that we only have numerical restriction on $I_i$\ (cf. (3.23))
and $d_i$\ can be small, thus we can not expect estimate of type
(3.17) to hold for all $\she\rdi$. Nevertheless, we have

\pro{Lemma 3.5}
There is a constant $C_1$\ depending only on $r$, $e_3>0$\ and
$\deg I\pri$, $I\pri\in\text{Pic}(X)$, such that for $r\pri\leq r$, we have
$$\mod\she^{r\pri,d\pri}_{e_s,I\pri}\leq \eta_X(r\pri,d\pri,I\pri)+C_1.
$$
\endpro
\def\rdp{^{r\pri,d\pri}_{e_3,I\pri}}

\proof
It suffices to show that there is a constant $C_1$\
such that for any $E\in\she\rdp$,
$$\dim\Ext^1(E,E)^0\leq 2r\pri d\pri-(r\pri-1)I^{\prime 2}+C_1.
$$
First, since $E$\ is $e_3$-stable and $e_3>0$, $\endo(E)$\
is $2re_3$-stable. Hence by lemma 1.8, both $h^0(\endo(E))$\ and
$h^0(\endo(E)\otimes K_X)$\ are bounded from above by a constant,
say $C_1$. By Serre duality, $\Ext^2(E,E)^0=H^0(\endo(E)\otimes K_X)$.
Therefore,
$$\align
\dim\Ext^1(E,E)^0&=2r\pri d\pri-(r\pri-1)I^{\prime 2}
-(r\pri-1)^2\xx(\OO_X)\\
&\qquad+\dim\Ext^0(E,E)^0+\dim\Ext^2(E,E)^0\\
&\leq 2r\pri d\pri-(r\pri-1)I^{\prime 2}+\Bl 2C_1-(r\pri-1)^2\xx(\OO_X)\Br.
\endalign
$$
This completes the proof of the lemma. \endpf

Returning to the proof of theorem 3.4, we need to estimate the term
$\dim\Ext^1(F_1,F_2)$\ in (3.25). First of all, by Riemann-Roch, for
$F_i\in\she\rdi$,
$$\align
\dim\Ext^1&(F_1,F_2)=\dim\Ext^0(F_1,F_2)+\dim\Ext^2(F_1,F_2)-\\
& -\Bigl({r_2\over 2}I_1^2+{r_1\over 2}I_2^2-
({r_1\over 2}I_2-{r_2\over 2}I_1)\cdot K_X-I_1\cdot I_2
+r_1r_2\xx(\OO_X)
-r_1 d_2-r_2 d_1\Bigr).\endalign
$$
Because $F_1$\ and $F_2$\ are $e_3$-stable,
$F_1^{\vee}\otimes F_2$\ and $F_2^{\vee}\otimes F_1$\  are $2re_3$-stable.
Also, the degree of $F_1^{\vee}\otimes F_2$\ and $F_2^{\vee}\otimes F_1$\
are bounded (from both sides) by
constants depending on $r$, $e$\ and $I\cdot H$. Thus, there is a constant
$C_2$\ depending on these parameters only so that
$$\dim\Ext^0(F_1,F_2),\quad \dim\Ext^2(F_1,F_2)\leq C_2.
$$
Therefore, for any admissible $(r_i, d_i, I_i)$,
$$\align
\mod\she&^{r_1,d_1}_{e_3,I_1}+\mod\she^{r_2,d_2}_{e_3,I_2}+
\sup\bigl\{\dim\Ext^1(F_1,F_2)\mid F_i\in\she\rdi\bigr\}\\
&\leq 2r_1d_1-(r_1-1)I_1^2
+ 2r_2d_2-(r_2-1)I_2^2+2C_1+\BBl-{r_2I_1^2-r_1 I_2^2 \over 2}+\\
&\qquad+
{r_1 I_2-{r_2}I_1\over 2}\cdot K_X+I_1\cdot I_2-r_1r_2\xx(\OO_X)
+r_1d_2+r_2d_1\BBr+2C_2\\
&\leq (2r-1)d+(1-r_2)d_1+(1-r_1)d_2-(r_1+{r_2\over 2}-1)I_1^2
-(r_2+{r_1\over 2}-1)I_2^2\\
&\quad +I_1\cdot I_2+({r_1\over 2}-{r_2\over 2}I_1)\cdot K_X-
r_1r_2\xx(\OO_X)+2C_1+2C_2.
\endalign
$$
Thanks to lemma 2.5, there is a constant $C_3\leq 0$\ depending on $r$,
$e$\ and
$I\cdot H$\ only such that $\displaystyle
d_i-{r_i-1\over 2r_i}I_i^2\geq C_3$.
Thus combined with $d=c_2(E)=I_1\cdot I_2+d_1+d_2$\ and $I_2=I-I_1$, the
right hand side of the above inequality is
$$\align
\leq& (2r-1)d-(r+{1-r_2\over 2r_1}-{3\over 2})I_1^2
-(r+{1-r_1\over 2r_2}-{3\over 2})I_2^2\tag 3.26\\
& -(2r-2)I_1\cdot I_2+({r_1\over 2}I_2
-{r_2\over 2}I_1)\cdot K_X+r^2|\xx(\OO_X)|+C_1+2C_2-rC_3\\
=&(2r-1)d+\BBl(1+{r_1-1\over 2r_2}+{r_2-1\over 2r_1})I_1^2+
{r_2+1-r_1\over 2r_2} I\cdot I_1-{r\over 2}I_1\cdot K_X\BBr+C_4.
\endalign
$$
Finally, because $|I_1\cdot H|\leq |I\cdot H|+e\sqrt{H^2}$, the
Hodge index theorem tells us that the sum of three middle terms in
the last line of (3.26) is bounded from above by a constant $C_5$.
Therefore, combined with (3.25), we have
$$\mod\bigl(\she\rd\setminus\she\rdz\bigr)\leq (2r-1)d+C\pri.   \eqno\qqed
$$

{\it Proof of (3.18)}:
We shall only consider the case where $e\geq 0$. The case $e<0$\ can
be proved similarly. First of all, by letting $e_1=e$\ and $e_2=0$\ in
theorem 3.4, we know that there is a constant $C_1$\ such that
$\mod\Bl\she\rd\setminus\she\rdz\Br
\leq (2r-1)d+C_1$. Then by choosing $N$\ large, we have $\mod\vecx\rdz=\erdi$\
and $(2r-1)d+C_1\leq\erdi$\ whenever $d\geq N$. Thus
$$\mod\vecx\rd\leq\max\big\{\mod\vecx\rdz,\mod\Bl\vecx\rd\setminus\vecx\rdz\Br
\big\}=\erdi.
$$
To prove (3.17), we will use the double
dual operation $\F$. Let
$$\F: \she\rd\lra \bigcup_{d\pri\leq d}\vecx^{r,d\pri}_{e,I}\tag 3.27
$$
be the map sending $E$\ to $E^{\vee\vee}$. Thanks to
lemma 3.2, we have
$$\mod\she\rd\leq\sup_{d\pri\leq d}\bigl\{\mod\vecx^{r,d\pri}_{e,I}
+(r+1)(d-d\pri)\bigr\}.
$$
Further, let $C_1\geq 0$\ be a constant such that
$$\mod\vecx^{r,d\pri}_{e,I}\leq\eta(r,d\pri,I)+C_1.
$$
Then, for $d\geq N+C_1$, $\mod\she\rd$\ is no more than either
$$\sup_{d\pri<N}\{\eta(r,d\pri,I)+C_1+(r+1)(d-d\pri)\}\leq \eta(r,N,I)+C_1\leq
\erdi
$$
or
$$\sup_{N\leq d\pri\leq d}\{\eta_X(r,d\pri,I)+(r+1)(d-d\pri)\}\leq\erdi.
$$
This establishes (3.17). (3.18) can be proved similarly
based on theorem 3.1.
We shall omit it.
\endpf

In light of the theorem 3.3, the proof of theorem 0.2 and 0.3 is
now quite easy. Recall that for the data $(r,d,I)$\ and sufficiently
large $n$, we can form the Grothendieck's
Quot-scheme $\bold{Quot}^{r,d}_{\rho,I}$\ of all
quotient sheaves $\OO_X^{\oplus
\rho}\to E$\ with $\rank E=r$, $\det E=I$, $c_2(E)=d$\ and $\rho=h^0(E(n))$.
If we let $\U\sub \bold{Quot}^{r,d}_{n,I}$\ be the open subset of all
semistable (with respect to $H$) quotient sheaves, then $\U$\ is
$SL(\rho,\CC)$\ invariant and the good quotient
$\U/\!/SL(\rho,\CC)$, which does exist,
is exactly the moduli scheme $\Mx$\ of rank $r$\ semistable sheaves of
$c_1=I$\ and $c_2=d$. Further, if we let $\U^s\sub\U$\ be the subset of
strictly stable sheaves, then $\pi\mh\U^s\to\pi(\U^s)\sub\Mx$\ is
a principal $SL(\rho,\CC)$-bundle. With this set up in mind, one
sees that in order to prove theorem 0.2, it suffices to classify the
singular locus of $\U$.
\def\bu{\bold{U}}

\pro{Proposition 3.6}
With the notation as before and for any constant $C$, there is a constant
$N$\ such that whenever $d\geq N$, then $\dim\U=\eta_X(r,d,I)+(\rho^2-1)$\
and the codimension
$$\codim\Bl\sing(\U),\U\Br\geq C.
$$
Further, when the codimension is at least 1, then $\U$\ locally is
a complete intersection and when the codimension is at least 2,
then $\U$\ is normal.
\endpro

\proof
Let $E\in\U$\ be any quotient sheaf, let
$q_2=h^0(\endo(E)\otimes K_X)$\ and let $q_1=\erdi+(\rho^2-1)+q_2$. Then
the argument in [L2, p8] demonstrates that the completion
of the local ring of
$\U$\ at $E$\ is of the form $k[[t_1,\cdots,t_{q_1}]]/J$, where $J$\ is
an ideal generated by at most $q_2$\ elements. In particular, for each
component $\bu\sub\U$, we always have
$$\dim\bu\geq \erdi+(\rho^2-1). \tag 3.28
$$
Next, by [At][Mu][Ma,p594], the singular locus $\sing(\U)$\ is exactly the set
of all quotient sheaves $E$\ with $\Ext^2(E,E)^0\ne0$. By theorem 3.3, for
any constant $C$, there is an $N$\ such that whenever $d\geq N$, the
set
$$\U\bigcap \she^{r,d}_{1,I}(K_X)=
\big\{E\in\U\mid h^0(\endo(E)\otimes K_X)\ne0\big\}
$$
obeys $\mod\Bl\U\cap\she^{r,d}_{1,I}(K_X)\Br\leq\erdi-C$.
Therefore,
$$
\dim \sing\Bl\U\Br
\leq\mod\BBl\U\bigcap \she^{r,d}_{1,I}(K_X)\BBr
+\dim SL(\rho) \leq \erdi+(\rho^2-1)-C.
$$
When $C\geq 1$, this inequality and (3.28) imply that $\U$\ has
purely dimension $\eta_X(r,d,I)+(\rho^2-1)$\ and
$\codim\Bl\sing(\U),\U)\Br\geq C$.
Because the completion of the local rings of $\U$\ are of
the form $k[[t_1,\cdots,t_{q_1}]]/J$\ with $J=(f_1,\cdots,f_{q_2})$. $\U$\
is a local complete intersection. $\U$\ will be normal if we further assume
$\codim\Bl\sing(\U),\U\Br\geq 2$.
\endpf

\pro{Corollary 3.7}
Let $X$\ be a smooth algebraic surface, $H$\ an ample divisor and $I$\
a line bundle on $X$. Let $r\geq 2$\ be an integer. Then for any constant
$C$, there is an $N$\ such that whenever $d\geq N$, then $\Mx$\ has pure
dimension $\erdi$\ and further, $\codim\Bl\sing(\mrdx),\mrdx\Br\geq C$.
\endpro

\proof
Since $\pi\mh\U^s\to\pi(\U^s)\sub\mrdx$\ is a principal bundle, the singular
locus $\sing(\mrdx)$\ is contained in
$$\pi\Bl\sing(\U^s)\Br\bigcap\pi\Bl\U\setminus\U^s\Br.
$$
By theorem 3.6, we know that for $d$\ large, we can arrange
$\codim\Bl\pi\Bl\sing(\U^s)\Br,\mrdx\Br\geq C$. Therefore, to prove the
corollary, we only need to find an upper bound of the dimension of
$\pi\Bl\U\setminus\U^s\Br$.

Let $E\in\U\setminus\U^s$. Then $E$\ admits a filtration
$0=E_0\sub E_1\sub \cdots\sub E_k=E$\
such that $F_i=E_i/E_{i-1}$\ are strictly stable. According to [Gi], $E$\
and $gr(E)=\oplus F_i$\ have the same image in $\mrdx$\ under $\pi$.
Thus $\dim\pi(\U\setminus\U^s)$\ can be bounded easily in terms of the
dimension of moduli of lower rank stable sheaves. Similar to the proof
of lemma 3.4, we can show that there is a constant $C_1$\ such that
$$\dim\pi(\U\setminus\U^s) \leq (2r-1)d+C_1.\tag 3.29
$$
(If we let $\she\rd$\ be the set introduced in \S 2, $\pi(\U\setminus\U^s)
\sub\she^{r,d}_{-1,I}\setminus\she^{r,d}_{\mu,I}$\ and then (3.29) follows
directly from lemma 3.4.) Thus for large $N$, we will have for $d\geq N$,
$\dim\pi(\U\setminus\U^s)\leq \erdi-C$. This completes the proof of the
corollary and the theorem 0.2.
\endpf

\pro{Corollary 3.8}
With the notation as before, then there exists $N$\ such that
whenever $d\geq N$, then
\roster
\item $\mrdx$\ is normal. Further, if $s\in\mrdx$\
is any closed point
corresponds to a stable sheaf, then $\mrdx$\ is a local
complete intersection at $s$. \hfil\break
\item The set of locally free $\mu$-stable sheaves $\Bl\mrdx\Br^{\text{vb}}\sub
\mrdx$\ is dense in $\mrdx$.
\endroster
\endpro

\proof
Let $N$\ be given by proposition 3.6 so that whenever $d\geq N$, $\U$\
has purely dimension $\eta_X(r,d,I)+(\rho^2-1)$\ and
$\codim\Bl\sing(\U),\U\Br\geq 2$. Then since $\U$\ is normal,
$\mrdx$\ must be normal and since $\U^s$\
is a local complete intersection, $\pi(\U^s)
\sub\mrdx$\ must be a local complete intersection. Here we have used the fact
that $\U\to\mrdx$\ is
a good quotient and $\U^s\to\pi(\U^s)$\ is a principal bundle.
The last statement can be proved easily similar to that of
theorem 3.1. We shall omit it here.
\endpf

The last subject we will study is the dependence of the moduli spaces on
the choice of the polarizations. We prove

\pro{Theorem 3.9}
For any choice $(r,I)$\ and polarizations $H_1$\ and $H_2$, there is a constant
$N$\ so that whenever $d\geq N$, then $\frak M^{r,d}_X(I,H_1)$\ and
$\frak M^{r,d}_X(I,H_2)$\ are birational to each other.
\endpro

\proof
Let $W\sub \frak M^{r,d}_X(I,H_1)$\ be the set of quotient sheaves
$E$\ such that $E$\ are not $H_2$-stable. Then every $E\in W$\ belongs to
the exact sequence
$$0\lra F_1\lra E\lra F_2\lra 0
$$
such that $\mu(F_1)\geq \mu(E)$\ (with respect to $H_2$).
Then by repeating the argument in lemma 3.4, we can find a constant $C_1$\
(depending on $H_1$\ and $H_2$) such that $\dim W\leq (2r-1)d+C_1$.
Therefore, by letting $N$\ large, we will have
$$\dim \frak M^{r,d}_X(I,H_1)=\dim\frak M^{r,d}_X(I,H_2)=\erdi
$$
and $\dim W\leq\erdi-1$\
provided $d\geq N$. Therefore, by the universality of the moduli
scheme, there is a morphism
$$\Phi:\frak M^{r,d}_X(I,H_1)\setminus W\lra \frak M^{r,d}_X(I,H_2)
$$
which is generically one-one and onto. Thus $\frak M^{r,d}_X(I,H_1)$\
is birational to $\frak M^{r,d}_X(I,H_2)$.
\endpf

\end